\documentclass[11pt]{article}
\usepackage{pdfsync}
\usepackage{hyperref}
\usepackage{array} 
\usepackage{amsfonts}
\usepackage{amsmath}
\usepackage{amssymb}
\usepackage{graphicx}
\usepackage{color}

\newtheorem{thm}{Theorem}[section]
\newtheorem{cor}[thm]{Corollary}
\newtheorem{lem}[thm]{Lemma}

\newtheorem{defn}[thm]{Definition}

\numberwithin{equation}{section}

\newcommand{\RR}{\mathbb{R}}
\newcommand{\NN}{\mathbb{N}}

\def\dist{\mathrm{dist}} 

\def\qed{\,\unskip\kern 6pt \penalty 500d
\raise -2pt\hbox{\vrule \vbox to8pt{\hrule width 6pt
\vfill\hrule}\vrule}\par}

\begin{document}

\title{\textbf{Travelling wave behaviour arising in nonlinear \\diffusion
problems posed in tubular domains} \\[3mm]}

\author{{\Large Alessandro Audrito\footnote{Dipartimento di Matematica ``Giuseppe Luigi Lagrange'' (DISMA), Politecnico di Torino, Italy.}  \;\, and \; Juan Luis V\'azquez\footnote{Departamento de Matem\'{a}ticas, Universidad Aut\'{o}noma de Madrid, Spain.}}
} 
\date{\vspace{-5ex} }

\maketitle

\begin{abstract}
For a fixed bounded domain $D \subset \RR^N$  we investigate the asymptotic behaviour for large times of solutions to the $p$-Laplacian diffusion equation posed in a tubular domain
\[
\partial_tu = \Delta_p u \quad \text{ in } D \times \RR, \quad t > 0
\]
with $p>2$, i.e., the  slow diffusion case, and homogeneous Dirichlet boundary conditions on the tube boundary. Passing to suitable re-scaled variables, we show the existence of a travelling wave solution  in logarithmic time \normalcolor that connects the level $u = 0$ and the unique nonnegative steady state associated to the re-scaled problem posed in a lower dimension, i.e. in $D\subset \RR^N$.

We then employ this special wave to show that a wide class of solutions converge to the universal stationary profile in the middle of the tube and at the same time they spread in both axial tube directions, miming the behaviour of the travelling wave (and its reflection) for large times.

The first main feature of our analysis is that wave fronts are constructed through a (nonstandard) combination of diffusion and absorbing boundary conditions, which gives rise to a sort of Fisher-KPP long-time behaviour. The second one is that the nonlinear diffusion term plays a crucial role in our analysis. Actually, in the linear diffusion framework $p=2$ solutions behave quite differently.
\end{abstract}


%
%
%
%
%
%
%
%
%
%
\section{Introduction}\label{SectionIntroduction}
We study the long-time behaviour of nonnegative solutions $u(x,t)$ to the initial value problem with $p$-Laplacian diffusion
\begin{equation}\label{eq:PLETUBULARDOMAIN}
\begin{cases}
\partial_tu = \Delta_p u \quad &\text{ in } \Omega \times 0,\infty) \\
u = 0 \quad &\text{ in } \partial\Omega \times (0,\infty) \\
u(x,0) = u_0(x) \quad &\text{ in } \Omega,
\end{cases}
\end{equation}
posed in a tubular domain $\Omega \subset \RR^{N+1}$ of the form
\begin{equation}\label{eq:TUBULARDOMAINASSUMPTION}
\Omega = D \times \RR,
\end{equation}
where $D \subset \RR^N$ is a bounded domain with smooth boundary.  Here, $\Delta_p u$ denotes the usual $p$-Laplace operator,  see motivations and details below, and we take $p>2$ (slow diffusion case).  We apply homogeneous Dirichlet boundary conditions. Such type of conditions are important for the type of results we obtain. The initial data are assumed to satisfy
\[
u_0 \geq 0 \quad \text{in } \Omega, \qquad u_0 \in \mathcal{C}_c(\Omega).
\]
We will use the notation $x=(y,z)$ with $y\in\RR$ and $z\in D$ for points $x\in \Omega$.

\paragraph{Travelling waves in tubes.} The main novelty of this paper lies in the peculiar form of the spatial domain $\Omega$, i.e., a tube. In the past years significant effort has been devoted to the study of evolution equations (linear and not) posed in bounded domains or in the whole space and the problem of describing the asymptotic behaviour of nonnegative solutions for large times is pretty well understood (see for instance \cite{Ag-Blan-Car:art,AronsonPeletier1985:art, KaminVazquez1988:art,S-V1:art,V2:book,Vazquez2004:art}). What we will describe below is in sharp contrast with what happens for solutions to \eqref{eq:PLETUBULARDOMAIN} when $\Omega \subset \RR^N$ is bounded and/or when $\Omega = \RR^N$. In the first case, nonnegative solutions converge to the unique nonnegative solution in separate variables form (cf. with \cite[Theorem 2.1]{S-V1:art}) while, in the second, they converge to a uniquely defined nonnegative solution in self-similar form (cf. with \cite[Theorem 2]{KaminVazquez1988:art}).

The present setting is less studied and it shows interesting differences both from the physical (combustion theory, fluid dynamics and so on \cite{C-D-D-S-V:art,Kal:survey,L-S-U:book}) and the mathematical viewpoint: we will show that nonnegative solutions to \eqref{eq:PLETUBULARDOMAIN} behave as a \emph{travelling wave} (TW for short) for large times, when computed at the correct \rm re-scaled variables (cf. with Theorem \ref{ASYMPTOTICBEHAVIOURTHEOREM}). In particular, the new time scale is logarithmic.
Let us point out that our investigation was motivated by the study of tubular propagation for the Porous Medium Equation done by V\'azquez in \cite{Vazquez2004:art,Vazquez2007:art}. We prove that the $p$-Laplacian setting also produces the logarithmic-time TW behaviour in tubes.

We recall that wave fronts are very common phenomena in reaction-diffusion equations and systems (see for instance the books \cite{GKbook, Smoller, Volpert3}, the classical papers \cite{Aro-Wein1:art,FifeMcLeod1977:art,Fisher:art,K-P-P:art}, while  \cite{DP-S:art,DePablo-Vazquez1:art,Eng-Gav-San:art} for the nonlinear diffusion framework). In previous works \cite{AA2017:art,AA-JLV:art,AA-JLV2:art}, the authors have studied the existence/non-existence of travelling fronts for reaction-diffusion equations with doubly nonlinear diffusion (slow and fast diffusion) and the asymptotic behaviour of more general solutions. Doubly nonlinear diffusion models include both porous medium and $p$-Laplacian models. Such results motivated us to investigate the existence/non-existence of TWs for similar nonlinear equations in different contexts.

In the present problem setting, the asymptotic behaviour for large times is in fact a sort of combination of two modes: a first mode has the separate variables form, and the second is a TW propagation. All solutions to \eqref{eq:PLETUBULARDOMAIN} (with nonnegative and compactly supported initial data) approach the special solution in separate variables form (corresponding to the domain $D \subset \RR^N$) on each compact sets of $\Omega$, while they show a travelling wave behaviour at both ends of the tube. \normalcolor The reason the TWs appear also in the study of a purely diffusive initial-value problem like \eqref{eq:PLETUBULARDOMAIN} is due to the shape of the domain. When $\Omega \subset \RR^{N+1}$ is a tube, the solutions spread along the longitudinal variable and the combination of the original diffusion with the reaction term (appearing in the re-scaled problem due to the change of variables), plus the loss of mass through the absorbing boundary $\partial \Omega$ gives rise to a travelling wave behaviour.
\paragraph{The $\boldsymbol{p}$-Laplacian.} We recall that the standard $p$-Laplacian operator is defined by the formula
\[
\Delta_p u := \nabla\cdot(|\nabla u|^{p-2}\nabla u),
\]
for smooth functions $u$. Here $\nabla$ and $\nabla \cdot$ denote the gradient and divergence operators (respectively) w.r.t. to the variable $x = (z,y) \in D\times\RR$. In some cases, we will use the symbols $\nabla_z$, $\nabla_z \cdot$ for the gradient and the divergence w.r.t. $z \in D$. The $p$-Laplacian has been widely studied as a standard model for elliptic and parabolic equations with degenerate-singular diffusion, and appears in a wide number of physical applications (see for instance \cite{C-D-D-S-V:art, DB:book, DiBenGianVes:book, Kal:survey, KaminVazquez1988:art, L-S-U:book, LindNotes, Urbano2008:art, V1:book} and huge number of references therein). In the present work we will take $p > 2$, corresponding to the degenerate diffusion range (cf. with \cite[Chapter 11]{V1:book}). The range $p < 2$ (singular diffusion) presents substantial differences and will not be studied here. We limit ourselves to present some general considerations at the end of the paper.
\paragraph{Free boundaries.} We finally recall that when $p > 2$ solutions $u = u(x,t)$ to \eqref{eq:PLETUBULARDOMAIN} (with compactly supported initial data) exhibit a \emph{free boundary}, i.e. the nonnegativity set of $u$ is a compact set of $\Omega$ for any $t > 0$. One of the main problems is thus to understand how it moves and which is its geometry, at least for large times. As a consequence of our main theorem, we will prove that the \emph{free boundary} of general solutions for large times is made of two finite sets described by two functions $y = s_{\pm}(z,t)$ that move towards both ends of the tube with logarithmic law:
\[
s_{\pm}(z,t) \sim \pm c_{\ast} \ln t, \quad \text{ uniformly in } z \in D, \quad \text{ for } t \sim +\infty,
\]
where $c_{\ast} = c_{\ast}(p,D) > 0$ is the speed of a critical travelling wave solution to problem \eqref{eq:REACTIONTRANSFORMATION} (cf.  Lemma \ref{Theorem:WAVECONSTRUCTION}). As claimed above, we thus obtain a nonstandard wave propagation behaviour, which strongly differs from both the case $\Omega \subset \RR^N$ and $\Omega = \RR^N$.


\paragraph{Structure of the paper.} The paper is divided in sections as follows.
The first three are introductory sections. In this section we have informally presented our work, motivating its interest and relating it with the existing literature. In Section \ref{SectionPreliminaries} we recall some known results concerning problem \eqref{eq:PLETUBULARDOMAIN} posed in bounded domains and two a priori estimates. In Section \ref{SectionMainResults} we introduce an equivalent transformed problem and we state our main result.

In Section \ref{Section:PreliminaryResults} we show two preliminary results concerning the asymptotic behaviour of solutions on compact sets of $\Omega$ and the existence of a family of wave solutions to
\begin{equation}\label{eq:EQUATIONREACTIONTRANSFORMATION}
\partial_{\tau} v = \Delta_p v + \frac{v}{p-2}   \quad \text{ in } \Omega \times (0,\infty),
\end{equation}
having zero Neumann derivative on the boundary $\partial \Omega$. Both of them play a fundamental role in the  construction of a special travelling wave solution to problem \eqref{eq:PLETUBULARDOMAIN} and in the proof of Theorem \ref{ASYMPTOTICBEHAVIOURTHEOREM}.

Section \ref{Section:ExistenceUniquenessTW} is the principal part of the work. We prove the existence (and we give information about the uniqueness) of a  wave solution to \eqref{eq:REACTIONTRANSFORMATION} connecting the solution $\Phi = \Phi(z)$ to \eqref{eq:STATIONARYPROBLEM} at $y = -\infty$ to the level $0$ at $y = +\infty$. This special solution will be employed as comparison tool in the proof of Theorem \ref{ASYMPTOTICBEHAVIOURTHEOREM}.

In Section \ref{Section:ProofMainTheorem} we present the proof of our main result. It is based on suitable comparison techniques with the wave solution constructed in Section \ref{Section:ExistenceUniquenessTW}.

Finally, in Section \ref{Section:Finalcomments} we present some open problems and future directions. The linear case $p = 2$ is also briefly discussed. Open directions are mentioned.

\section{Preliminaries}\label{SectionPreliminaries}
We recall needed facts on the  behaviour of solutions in bounded domains as a basic preliminary of what follows. We introduce weak solutions, and we prove basic a priori estimates.
\paragraph{Some properties on bounded domains.} As we have said above, we recall now a crucial result concerning the long time behaviour of solutions to \eqref{eq:PLETUBULARDOMAIN} in the case in which the spatial domain is a bounded subset of $\RR^N$. So, let $D \subset \RR^N$ be a bounded domain with smooth boundary and consider the problem
\begin{equation}\label{eq:PLEBOUNDEDDOMAIN}
\begin{cases}
\partial_tu = \Delta_p u \quad &\text{ in } D \times (0,\infty) \\
u = 0 \quad &\text{ in } \partial D \times (0,\infty) \\
u(x,0) = u_0(x) \quad &\text{ in } D,
\end{cases}
\end{equation}
where $0 \leq u_0 \in L^1(D)$ is a nontrivial initial datum. In \cite{S-V1:art}, Stan and V\'azquez proved the following remarkable asymptotic behaviour theorem that we report here adapting it to our notations (note that they worked in the quite more general context of the ``doubly nonlinear diffusion''). For Porous Medium diffusion we quote the work of Aronson and Peletier \cite{AronsonPeletier1985:art}.
\begin{thm}(\cite[Theorem 2.1]{S-V1:art})\label{ASYMPTOTICSINBOUNDEDDOMAINS}
Let $u = u(z,t)$  be a nonnegative weak solution to problem \eqref{eq:PLEBOUNDEDDOMAIN}. Then
\[
\lim_{t \to \infty} t^{\frac{1}{p-2}}\|u(\cdot,t) - U(\cdot,t)\|_{L^{\infty}(D)} = 0,
\]
where $U = U(z,t)$ is the unique solution to
\begin{equation}\label{eq:PLEBOUNDEDDOMAINSEPVARIABLE}
\begin{cases}
\partial_tU = \Delta_p U \quad &\text{ in } D \times (0,\infty) \\
U = 0 \quad &\text{ in } \partial D \times (0,\infty),
\end{cases}
\end{equation}
in separate variables form
\[
U(z,t) = t^{-\frac{1}{p-2}} \,\Phi(z),
\]
where $\Phi = \Phi(z)$ is the unique nonnegative and nontrivial weak solution to the stationary problem
\begin{equation}\label{eq:STATIONARYPROBLEM}
\begin{cases}
-(p-2)\Delta_p \Phi = \Phi \quad &\text{ in } D \\
\Phi = 0              \quad &\text{ in } \partial D.
\end{cases}
\end{equation}
\end{thm}

The proof of the above theorem is based on a clever change of variables, together with some monotonicity results and a priori estimates (we will employ some of these methods later).
\paragraph{Weak solutions and two fundamental a priori estimates.} This paragraph is devoted to the introduction of the basic concepts and ideas of paper. We report below the definition of weak solutions to problem \eqref{eq:PLETUBULARDOMAIN}. Weak solutions exist and are unique (see for instance \cite{DBen-Her2:art}). Note that thanks to well-known regularity results (\cite{DB:book,DB-Fried:art,DiBenGianVes:book,Urbano2008:art}), we can assume that solutions to \eqref{eq:PLETUBULARDOMAIN} are $C^{1,\alpha}(\Omega\times[s,\infty))$ for all $s > 0$ and some $\alpha \in (0,1)$, and furthermore, the equation is satisfied for a.e. $(x,t) \in \Omega\times(0,\infty)$.

\begin{defn}
A nonnegative weak solution to problem \eqref{eq:PLETUBULARDOMAIN} is a nonnegative function $u = u(x,t)$ satisfying the following properties:
\begin{itemize}
  \item $u \in \mathcal{C}([0,\infty):L_{loc}^1(\Omega))$ and $u(t) \to u_0$ in $L^1(\Omega)$ as $t \to 0$.
  \item $u \in L_{loc}^1(0,\infty:W_{0,loc}^{1,p}(\Omega))$.
  \item The identity
\[
\int_0^{\infty}\int_{\Omega}\left\{|\nabla u|^{p-2}\nabla u \cdot \nabla \eta - u \partial_t\eta \right\} dxdt = 0
\]
holds for every $\eta \in \mathcal{C}_c^{\infty}(\Omega\times(0,\infty))$.
\end{itemize}
\end{defn}

The work relies on the validity two important \emph{a priori} estimates. The first one, is the so called ``universal estimate'':
\begin{equation}\label{eq:UNIVERSALESTIMATE}
0 \leq u(x,t) \leq C \, t^{-\frac{1}{p-2}} \quad \text{ in } \Omega \times (0,\infty),
\end{equation}
where $C = C(N,p) > 0$, while the second is the well-known B\'enilan-Crandall type estimate
\begin{equation}\label{eq:BENILANCRANDALLESTIMATE}
\partial_t u \geq - \frac{u}{(p-2)t} \quad \text{ in } \Omega \times (0,\infty),
\end{equation}
satisfied by nonnegative weak solutions $u = u(x,t)$ to \eqref{eq:PLETUBULARDOMAIN} in the sense of distributions. Note that \eqref{eq:BENILANCRANDALLESTIMATE} implies that the function $t \to t^{1/(p-2)} u(x,t)$ is nondecreasing. Let us shortly review their proofs.

To prove \eqref{eq:UNIVERSALESTIMATE} we proceed in two steps. First of all, we compare the solution $u = u(x,t)$ to \eqref{eq:PLETUBULARDOMAIN} with the solution $\widetilde{u} = \widetilde{u}(x,t)$ to the Cauchy problem posed in the all space
\[
\begin{cases}
\partial_t \widetilde{u} = \Delta_p \widetilde{u}  \quad &\text{ in } \RR^{N+1}\times(0,\infty) \\
\widetilde{u}(x,0) = \widetilde{u}_0(x)                        \quad &\text{ in } \RR^{N+1},
\end{cases}
\]
where $\widetilde{u}_0 = u_0$ in $\Omega$, while $\widetilde{u}_0 = 0$ in $\RR^{N+1}\setminus\Omega$. Since $\Omega \subset \RR^{N+1}$, it turns out
\[
0 \leq u(x,t) \leq \widetilde{u}(x,t) \leq C \|\widetilde{u}_0\|_{L^1(\RR^{N+1})}^{\frac{p\alpha}{N+1}} t^{-\alpha}, \qquad \alpha = \frac{N+1}{(N+1)(p-2) + p},
\]
for some constant $C = C(N,p)$. Note that the second inequality of the above chain is justified by the comparison principle for ``domain variations'' (cf. with Proposition 6.9 of \cite{V2:book} for the Porous Medium setting), while the third one, by the ``smoothing effect'' estimate (cf. with Theorem 11.3 of \cite{V1:book}). In particular, it follows that $u(\cdot,t)$ is bounded for all $t > 0$. In the second step, we improve the above estimate by comparing the solution $u = u(x,t)$ \eqref{eq:PLETUBULARDOMAIN} with the solution $\widehat{u} = \widehat{u}(z,t)$ to problem \eqref{eq:PLEBOUNDEDDOMAIN} (starting at $\varepsilon$):
\[
\begin{cases}
\partial_t\widehat{u} = \Delta_p \widehat{u} \quad &\text{ in } D \times (\varepsilon,\infty) \\
\widehat{u} = 0 \quad &\text{ in } \partial D \times (\varepsilon,\infty) \\
\widehat{u}(z,\varepsilon) = \widehat{u}_{\varepsilon} \quad &\text{ in } D,
\end{cases}
\]
where $\widehat{u}_{\varepsilon} \equiv \|u(\cdot,\varepsilon)\|_{L^{\infty}(\Omega)}$, for some $\varepsilon > 0$ fixed (note that we are crucially exploiting the fact the $u = u(x,t)$ is bounded for any $t > 0$). Consequently, since $\widehat{u} = \widehat{u}(z,y,t)$ can be seen as a solution to \eqref{eq:PLETUBULARDOMAIN} which is constant w.r.t. $y \in \RR$ and $u(z,y,\varepsilon) \leq \widehat{u}_{\varepsilon}$, it follows by comparison $u(z,y,t+\varepsilon) \leq \widehat{u}(z,y,t+\varepsilon)$ for all $t \geq 0$. Consequently, taking $t=0$ and using the arbitrariness of $\varepsilon > 0$, we get
\[
u(x,t) \leq \widehat{u}(x,t) \leq \Phi(z) t^{-\frac{1}{p-2}} \leq C t^{-\frac{1}{p-2}}, \quad \text{ in } \Omega \times (0,\infty),
\]
where the function $\Phi = \Phi(z)$ in the second inequality is defined in Theorem \ref{ASYMPTOTICSINBOUNDEDDOMAINS} (cf. with \cite{S-V1:art} and Chapter 5 of \cite{V2:book} for the Porous Medium setting). This concludes the proof of \eqref{eq:UNIVERSALESTIMATE}.

For what concerns \eqref{eq:BENILANCRANDALLESTIMATE}, we give a nice proof based on scaling, valid for a smooth solution $u = u(x,t)$ to problem \eqref{eq:PLETUBULARDOMAIN} (the proof for weak solutions can be done with a suitable approximation technique). Let us consider the family of functions
\[
u_{\lambda}(x,t) = \lambda u(x,\lambda^{p-2}t),
\]
for any value of the parameter $\lambda \geq 1$. It is easily seen that $u_{\lambda} = u_{\lambda}(x,t)$ is still a solution to problem \eqref{eq:PLETUBULARDOMAIN} with $u_{\lambda}(x,0) = \lambda u(x,0) \geq u(x,0)$, since $\lambda \geq 1$. We thus obtain by comparison $u(x,t) \leq u_{\lambda}(x,t)$ in $\Omega\times(0,\infty)$, for all $\lambda \geq 1$. Hence,
\[
0 \leq \frac{\partial u_{\lambda}}{\partial \lambda} \, \Big|_{\lambda = 1} = u + t(p-2)\partial_tu \quad \text{ in } \Omega \times (0,\infty),
\]
from which we immediately deduce \eqref{eq:BENILANCRANDALLESTIMATE}.
\section{Transformed problem and main results}\label{SectionMainResults}
\paragraph{New Problem.} We now introduce one of the main ideas of the work. It is based on the change of variables (that we refer as scaling of re-normalization)
\begin{equation}\label{eq:RENORMALIZATIONFORMULA}
v(x,\tau) = t^{\frac{1}{p-2}}u(x,t), \qquad \tau = \ln t,
\end{equation}
which transforms problem \eqref{eq:PLETUBULARDOMAIN} into the reaction-diffusion problem with homogeneous Dirichlet boundary conditions
\begin{equation}\label{eq:REACTIONTRANSFORMATION}
\begin{cases}
\partial_{\tau} v = \Delta_p v + \frac{v}{p-2}   \quad &\text{ in } \Omega \times \RR \\
v = 0 \quad &\text{ in } \partial\Omega \times \RR,
\end{cases}
\end{equation}
in the sense that $v \in L_{loc}^p(\RR:W_{0,loc}^{1,p}(\Omega))$ and, furthermore, the identity
\[
\int_{-\infty}^{\infty}\int_{\Omega}\left\{|\nabla v|^{p-2}\nabla v \cdot \nabla \eta - v \left( \frac{\eta}{p-2} + \partial_{\tau}\eta\right) \right\} dx d\tau = 0,
\]
holds for every $\eta \in \mathcal{C}_c^{\infty}(\Omega\times\RR)$. It thus follows that the equations in \eqref{eq:PLETUBULARDOMAIN} and \eqref{eq:REACTIONTRANSFORMATION} (together with the homogeneous Dirichlet boundary conditions) are completely equivalent and linked by relation \eqref{eq:RENORMALIZATIONFORMULA}. Nevertheless, we have to stress that the function $v = v(x,\tau)$ is an \emph{eternal solution} (i.e. defined for any time $\tau \in \RR$), making impossible a precise formulation of the initial condition. This inconvenient can be easily overcome by noting that we can modify transformation \eqref{eq:RENORMALIZATIONFORMULA} to
\begin{equation}\label{eq:SECONDRENORMALIZATIONFORMULA}
v(x,\tau) = (t+1)^{\frac{1}{p-2}}u(x,t), \qquad \tau = \ln (t+1),
\end{equation}
to have $\tau \in [0,\infty)$ and $v(x,0) = u_0(x)$ (i.e. $v \in \mathcal{C}([0,\infty):L_{loc}^1(\Omega))$ and $v(\tau) \to u_0$ in $L^1(\Omega)$ as $\tau \to 0$). Since we are only concerned with the asymptotic behaviour for $t \to +\infty$ (i.e. $\tau \to +\infty$), we have chosen to work with the transformation in \eqref{eq:RENORMALIZATIONFORMULA} which has the advantage to keep the notations as simpler as possible.

As a first important byproduct, it turns out that both universal and B\'enilan-Crandall estimates \eqref{eq:UNIVERSALESTIMATE} and \eqref{eq:BENILANCRANDALLESTIMATE} are transformed into
\begin{equation}\label{eq:UNIVERSALESTIMATE1}
0 \leq v \leq C \quad \text{ in } \Omega \times \RR,
\end{equation}
where $C = C(N,p) > 0$ is a new constant, and
\begin{equation}\label{eq:BENILANCRANDALLESTIMATE1}
\partial_{\tau} v \geq 0 \quad \text{ in } \Omega \times \RR,
\end{equation}
respectively. This new way of looking at \eqref{eq:UNIVERSALESTIMATE} and \eqref{eq:BENILANCRANDALLESTIMATE} is crucial in the asymptotic behaviour analysis. Indeed, from \eqref{eq:UNIVERSALESTIMATE1} and \eqref{eq:BENILANCRANDALLESTIMATE1} we immediately deduce the existence of a bounded point-wise limit $\overline{\Phi} = \overline{\Phi}(x)$:
\[
v(\cdot,\tau) \to \overline{\Phi}(\cdot) \quad \text{ in } \Omega \; \text{ as } \tau \to +\infty.
\]
The above convergence is the first step in our long time behaviour analysis and completely motivates the re-normalization \eqref{eq:RENORMALIZATIONFORMULA}. We will thus prevalently work with the transformed solution $v = v(x,\tau)$ and then recover the information on $u = u(x,t)$ through \eqref{eq:RENORMALIZATIONFORMULA}-\eqref{eq:SECONDRENORMALIZATIONFORMULA}.
\paragraph{Main result.} We have reported above some known facts and ideas which are fundamental ingredients of this work. We now state our main theorem which describes in a precise and quantified way the long time behaviour of solutions $u = u(z,y,t)$ to problem \eqref{eq:PLETUBULARDOMAIN}. For clarity we state it for solutions $v = v(z,y,\tau)$ to problem \eqref{eq:REACTIONTRANSFORMATION}, obtaining the information on $u = u(z,y,t)$ through the change of variables \eqref{eq:RENORMALIZATIONFORMULA}.
\begin{thm}\label{ASYMPTOTICBEHAVIOURTHEOREM}
There exists a unique number $c_{\ast} = c_{\ast}(p,D) > 0$ such that problem \eqref{eq:REACTIONTRANSFORMATION} admits a continuous nonnegative weak solution with wave form
\[
v(z,y,t) = \varphi(z,y - c_{\ast}\tau),
\]
which is unique (in the sense of Lemma \ref{Lemma:Uniqueness}) up to shifts along the longitudinal variable, non-increasing w.r.t. $y \in \RR$ and satisfies
\begin{equation}\label{eq:LIMITCONDTWPROFILES0}
\lim_{\xi \to -\infty} \varphi(z,\xi) = \Phi(z), \qquad \lim_{\xi \to +\infty} \varphi(z,\xi) = 0, \quad z \in D,
\end{equation}
where $\Phi = \Phi(z)$ is the unique nonnegative weak solution to problem \eqref{eq:STATIONARYPROBLEM} (cf. with Theorem \ref{ASYMPTOTICSINBOUNDEDDOMAINS}).

Furthermore, any solution $v = v(z,y,\tau)$ to \eqref{eq:REACTIONTRANSFORMATION} with nontrivial and nonnegative compactly supported initial datum $v_0$ satisfies the following assertions:

(i) For all $0 < c < c_{\ast}$,
\[
\limsup_{\tau \to +\infty} \sup_{z \in D, \; |y| \leq c\tau} |v(z,y,\tau) - \Phi(z)| = 0.
\]

(ii) For all $c > c_{\ast}$, there exists a time $\tau_c > 1$ such that
\[
v(z,y,\tau) = 0 \quad \text{ in } \{z \in D\}\times\{|y| \geq c\tau\}, \quad \text{for all } \tau \geq \tau_c.
\]
Finally, there exists a waiting time $T > 0$ (depending only on the data and $D$) such that for all $\tau > T$, the free boundary of $v = v(z,y,\tau)$ is made by two disjoint sets
\[
S_v^{\pm}(\tau) := \{(z,S_v^{\pm}(z,\tau)): z \in D \},
\]
where
\[
S_v^+(z,\tau) := \inf\{y > 0: v(z,y,\tau) = 0 \}, \qquad S_v^-(z,\tau) := \sup\{y < 0: v(z,y,\tau) = 0 \}
\]
are locally smooth functions for any $z \in D$, $\tau > T$ and
\[
\lim_{\tau \to +\infty} \frac{S_v^+(z,\tau)}{\tau} = - \lim_{\tau \to +\infty} \frac{S_v^-(z,\tau)}{\tau} = c_{\ast}.
\]
\end{thm}
The above statement needs some comments. First of all, we notice that, passing to the initial variables, we have proven the existence of a special solution to \eqref{eq:PLETUBULARDOMAIN} which is a re-scaled and re-normalized TW:
\[
u(z,y,t) = t^{-\frac{1}{p-2}} \varphi(z,y - c_{\ast}\ln t),
\]
and $\varphi = \varphi(z,\xi)$ satisfies \eqref{eq:LIMITCONDTWPROFILES0}. Similar, part (i) and (ii) can be reformulated as
\[
\limsup_{t \to +\infty} \; t^{\frac{1}{p-2}} \sup_{z \in D, \; |y| \leq c\ln t} |u(z,y,t) - U(z,t)| = 0,
\]
for any fixed $0 < c < c_{\ast}$, where $U = U(z,t)$ is the unique nonnegative solution to \eqref{eq:PLEBOUNDEDDOMAINSEPVARIABLE} in separate variables form (cf. with Theorem \ref{ASYMPTOTICSINBOUNDEDDOMAINS}), and
\[
u(z,y,t) = 0 \quad \text{ in } \{z \in D\}\times\{|y| \geq c\ln t\}, \quad \text{for all } t \geq t_c,
\]
for any fixed $c > c_{\ast}$ and some $t_c > 0$ large enough (depending on $c$). Some simulations of such behaviour are displayed in Figure \ref{fig:SimulationsTwoBumps} below.

\begin{figure}[!ht]
  \centering
  \includegraphics[scale = 0.8]{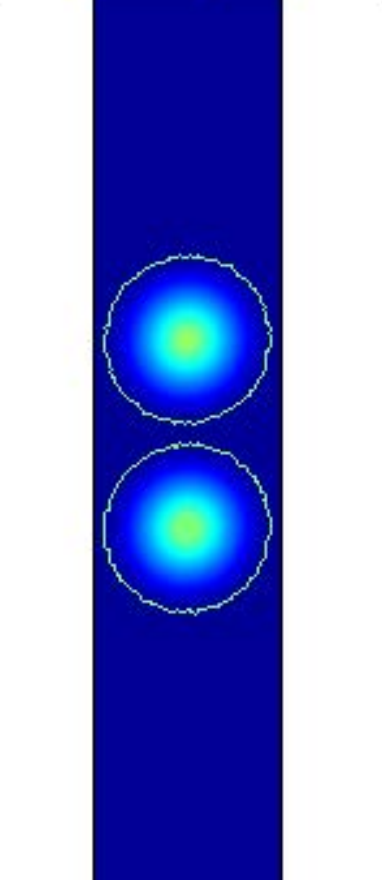} \quad
  \includegraphics[scale = 0.8]{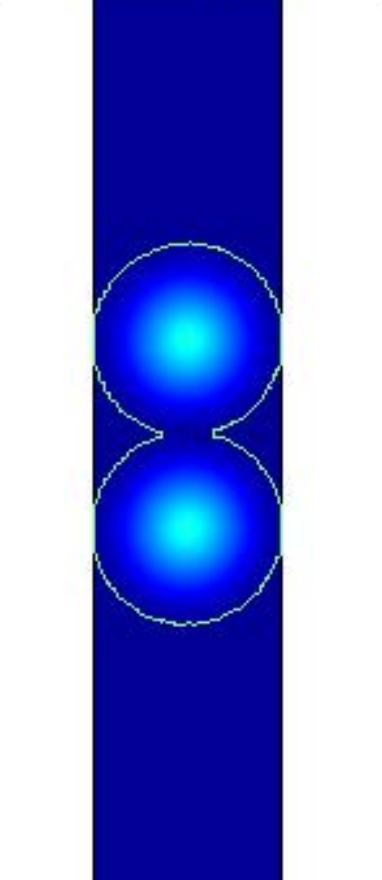} \quad
  \includegraphics[scale = 0.8]{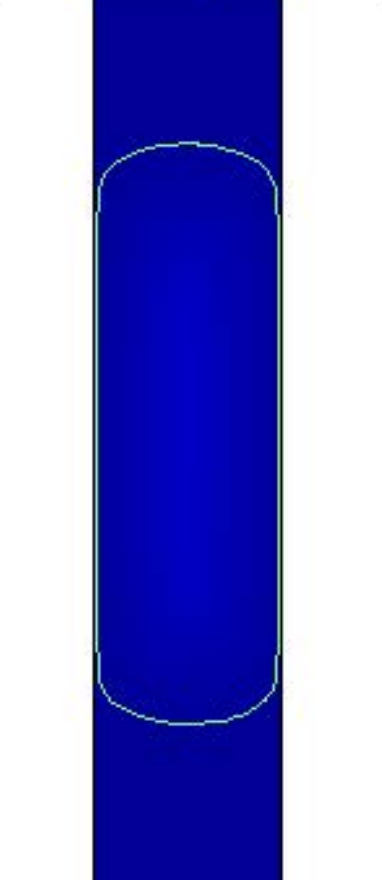}
  \caption{Three stages of the formation of the double travelling wave behaviour. The exponent is $p=4$, the space dimension is $1+1$, the initial data were taken to be nearly constant in two disjoint circles and zero outside, and time is logarithmic (simulation software: COMSOL).}\label{fig:SimulationsTwoBumps}
\end{figure}

We have thus stated in a precise way, what was anticipated in the introduction. First, the solution in separate variables form $U = U(z,t)$ (cf. with Theorem \ref{ASYMPTOTICSINBOUNDEDDOMAINS}) is stable and is the uniform limit for large times of solutions $u = u(z,y,t)$ to \eqref{eq:PLETUBULARDOMAIN} if we restrict ourselves to compact sets of $\Omega$. Quoting V\'azquez \cite{Vazquez2007:art}: ``the intermediate asymptotic behaviour of our problem in a tube forgets the longitudinal variable and decays in first approximation like the transversal problem in one dimension less''. On the other hand, there is a propagation phenomenon. Part (i) explains that the convergence to the special solution $U = U(z,t)$ takes place not only on compact sets of $\Omega$, but also in all subsets of the type
\[
D \times \{|y| \leq c \ln t\}
\]
for large times, where $c$ is arbitrarily fixed between $0$ and $c_{\ast}$. Part (ii) assures that the front of the solution (and its \emph{free boundary}) cannot move faster than $c\ln t$ when $t$ is large and $c > c_{\ast}$ is fixed. As a consequence of part (i) and (ii) we obtain the description of the wave behaviour of solutions to \eqref{eq:PLETUBULARDOMAIN} for large times.

As we have anticipated in the introduction, we can give an interesting interpretation to Theorem \ref{ASYMPTOTICBEHAVIOURTHEOREM}, based on a comparison with some results on reaction-diffusion equations we have recently obtained. In a recent paper \cite{AA-JLV:art} (see also \cite{AA2017:art}) we constructed a family of TWs $v(x,\tau) = \varphi(n \cdot x - c\tau)$ to the Fisher-KPP equation
\begin{equation}\label{eq:KPPpLaplacian}
\partial_{\tau} v + \Delta_p v = f(v) \quad \text{ in } \RR^N\times(0,\infty),
\end{equation}
where $p > 2$, $c > 0$, $n$ is a unit vector and $f$ is modeled on the logistic reaction term $f(v) = v (1 - v)$, and we employ them to study the long-time behaviour of solutions with initial data $0 \leq v_0 \leq 1$, with compact support. We showed the existence of a minimal speed of propagation $c_{\ast} > 0$ (and a corresponding TW with free boundary) such that any general solution satisfies
\[
\limsup_{\tau \to +\infty} \sup_{|x| \leq c\tau} |v(x,\tau) - 1| = 0,
\]
for any fixed $0 < c < c_{\ast}$, while
\[
v(x,\tau) = 0 \quad \text{ in } \{|x| \geq c\tau\}, \quad \text{for all } \tau \geq \tau_c,
\]
for any $c > c_{\ast}$ and some $\tau_c > 0$ (cf. with \cite[Theorem 2.6]{AA-JLV:art} ). This means that general solutions to \eqref{eq:KPPpLaplacian} (with compactly supported initial data $0 \leq v_0 \leq 1$) spread through the space with constant speed of propagation $c_{\ast} > 0$ for large times, converging to the (stable) steady state $v = 1$. Something similar happens in the present framework, where general solutions $v = v(z,y,\tau)$ to \eqref{eq:REACTIONTRANSFORMATION} converge to the unique nonnegative weak solution $\Phi = \Phi(z)$ to problem \eqref{eq:STATIONARYPROBLEM} with rate $c_{\ast} = c_{\ast}(p,D) > 0$. The main fact here is that the wave behaviour is generated in a different way: in the reaction-diffusion equations framework it is the byproduct of the combination of diffusion and reaction, while in this setting it comes from the interaction between diffusion, reaction and \emph{absorbing boundary}. Indeed, we notice that the reaction term $v$ in the r.h.s. of \eqref{eq:REACTIONTRANSFORMATION} is nothing more than the linear part of $f(v) = v(1 - v)$ and the role played by the term $-v^2$ in the Fisher-KPP setting is now played by the homogeneous Dirichlet conditions at the lateral boundary.  Let us finally stress that the existence of this TW solution is due to the presence of \emph{slow nonlinear diffusion}. Actually, the paper shows that in this setting the influence of the Dirichlet boundary conditions is much stronger than in the case of linear diffusion (an interesting case which is briefly discussed in  Section \ref{Section:Finalcomments}). \normalcolor
%
%
%
%
%
%
%
%
%
%
%
\section{Two preliminary results}\label{Section:PreliminaryResults}
As explained above, we begin our investigation with two preliminary results. The first one is a PDE lemma concerning the asymptotic behaviour of solutions to \eqref{eq:PLETUBULARDOMAIN} on compact sets of $\Omega$, while the second focuses on the description of a class of wave solutions to \eqref{eq:EQUATIONREACTIONTRANSFORMATION} with Neumann boundary conditions.

\paragraph{I. Convergence in compact sets of $\boldsymbol{\Omega}$.} The following lemma states that the long-time behaviour of nonnegative weak solutions to problem \eqref{eq:PLETUBULARDOMAIN} is independent of the $y$ variable, if we restrict our analysis to compact sets of $\Omega = D\times\RR$. Nonnegative weak solutions converge to the special solution in separate variables form of the Cauchy-Dirichlet problem posed in the bounded domain $D \subset \RR^N$.
\begin{lem}\label{LEMMACONVERGENCEINBOUDEDDOMAINS}
Let $u = u(z,y,t)$ be a nonnegative weak solution to problem \eqref{eq:PLETUBULARDOMAIN} with nontrivial and nonnegative initial data $u_0 \in L^1(\Omega)$. Then
\[
\lim_{t \to \infty} t^{\frac{1}{p-2}}\|u(z,y,t) - U(z,t)\|_{L^{\infty}(K)} = 0,
\]
for any compact set $K \subset \Omega = D\times\RR$, where $U = U(z,t)$ is the unique solution to \eqref{eq:PLEBOUNDEDDOMAINSEPVARIABLE} in self-similar form
\[
U(z,t) = t^{-\frac{1}{p-2}} \,\Phi(z),
\]
where $\Phi = \Phi(z)$ is the unique nonnegative weak solution to the stationary problem \eqref{eq:STATIONARYPROBLEM} (cf. with Theorem \ref{ASYMPTOTICSINBOUNDEDDOMAINS}).
\end{lem}
\paragraph{Proof.} Let us fix a compact set $K \subset \Omega$. Following the proof of the Porous Medium case (see Theorem 3.2 of \cite{Vazquez2004:art}), we proceed in some few steps.

\emph{Step 1: Comparison from below.} We compare the solution $u = u(x,t)$ to problem \eqref{eq:PLETUBULARDOMAIN} with a sequence of functions $u_j = u_j(x,t)$ being weak solutions to
\[
\begin{cases}
\partial_tu_j = \Delta_p u_j \quad &\text{ in } \Omega_j \times (0,\infty) \\
u_j = 0 \quad &\text{ in } \partial\Omega_j \times (0,\infty) \\
u_j(x,0) = u_{0,j}(x) \quad &\text{ in } \Omega_j,
\end{cases}
\]
where $u_{0,j} = u_0$ in $\Omega_j$ and zero outside, whilst
\[
\Omega_j := D\times(-j,j) = \{(z,y): z \in D, \; |y| < j \},
\]
for any integer $j \geq 1$. By comparison it follows that $u_j = u_j(x,t)$ is a nondecreasing sequence in $j$ with
\[
u_j(x,t) \leq u(x,t) \quad \text{ in } \Omega_j\times(0,\infty),
\]
for all $j \geq 1$ and, furthermore, by Theorem \ref{ASYMPTOTICSINBOUNDEDDOMAINS} it follows
\[
\lim_{t \to \infty} \|t^{\frac{1}{p-2}} u_j(\cdot,t) - \Phi_j(\cdot)\|_{L^{\infty}(\Omega_j)} = 0,
\]
for any fixed $j \geq 1$, where $\Phi_j = \Phi_j(x)$ is the unique nonnegative weak solution to the stationary problem \eqref{eq:STATIONARYPROBLEM}:
\begin{equation}\label{eq:APPROXSEPVARSOLLEMMACONVONCOMPACT}
\begin{cases}
-(p-2)\Delta_p \Phi_j = \Phi_j \quad &\text{ in } \Omega_j \\
\Phi_j = 0              \quad &\text{ in } \partial \Omega_j,
\end{cases}
\end{equation}
i.e., $\Phi_j \in W_0^{1,p}(\Omega_j)$ and the identity
\begin{equation}\label{eq:EQUATIONFORPHIJ}
(p-2) \int_{\Omega_j} |\nabla \Phi_j|^{p-2} \nabla\Phi_j \cdot \nabla \phi \,dx = \int_{\Omega_j} \Phi_j \phi \,dx
\end{equation}
holds for every $\phi \in \mathcal{C}_c^{\infty}(\Omega_j)$. We notice that the above convergence is both uniform and monotone, and each $\Phi_j \in \mathcal{C}^{1,\alpha}(\Omega_j)$ from standard regularity theory (cf. for instance with \cite{DiBenedetto1983:art,Uhlenbeck1977:art,Uralceva1968:art}).

\emph{Step 2: Re-scaled orbit and convergence.} Now, through the change of variables  \eqref{eq:RENORMALIZATIONFORMULA}, we consider the re-normalized version of the solution $v = v(x,\tau)$, which, in view of \eqref{eq:UNIVERSALESTIMATE} and \eqref{eq:BENILANCRANDALLESTIMATE}, satisfies (cf. with the introduction)
\begin{equation}\label{eq:LIMITOFRENORMALIZEDVERSION}
v(\cdot,\tau) \to \overline{\Phi}(\cdot) \quad \text{ in } \Omega,
\end{equation}
with point-wise monotone convergence as $\tau \to +\infty$. The nondecreasing monotonicity implies that the limit $\overline{\Phi} = \overline{\Phi}(x)$ is not trivial and, secondly, that the convergence \eqref{eq:LIMITOFRENORMALIZEDVERSION} takes place in all $L^p(\Omega)$, for $1 \leq p < \infty$ (it sufficient to apply Beppo Levi monotone convergence theorem). Finally, note that from \emph{Step 1} and the proof of \eqref{eq:UNIVERSALESTIMATE}, we get also
\begin{equation}\label{eq:BOUNDONPHIJANDPHI}
\Phi_j(x) \leq \overline{\Phi}(x) \leq \Phi(z) \quad \text{ in } \Omega_j,
\end{equation}
for all $j \geq 1$.

\emph{Step 3: Identification of the limit.} Since the sequence $u_j = u_j(x,t)$ (introduced in \emph{Step 1}) is monotone nondecreasing in $j$ by construction and $t^{1/(p-2)}u_j(x,t)$ is nondecreasing in $t$ by \eqref{eq:BENILANCRANDALLESTIMATE}, we have that the sequence $\Phi_j = \Phi_j(x)$ is monotone nondecreasing in $j$ and, in view of \eqref{eq:BOUNDONPHIJANDPHI}, it follows
\[
\Phi_j \to \underline{\Phi} \quad \text{ in } \Omega
\]
as $j \to +\infty$ with point-wise convergence, for some limit function $\underline{\Phi} = \underline{\Phi}(x)$ satisfying
\begin{equation}\label{eq:INEQUALITYOFLIMITSFUNCTIONPHI}
\underline{\Phi}(x) \leq \overline{\Phi}(x) \leq \Phi(z) \quad \text{ in } \Omega.
\end{equation}
The monotonicity of $j \to \Phi_j$ implies that the above convergence holds in $L^p(\Omega)$, $1 \leq p < \infty$, and thus for any $\phi \in \mathcal{C}_c^{\infty}(\Omega)$,
\[
\int_{\Omega_j} \Phi_j \phi \,dx \to \int_{\Omega} \underline{\Phi} \phi \,dx, \quad \text{ as } j \to +\infty.
\]
Furthermore, $\{\Phi_j\}_{j\in\NN}$ is uniformly bounded in $L^{\infty}(\Omega)$ (this easily follows from \eqref{eq:BOUNDONPHIJANDPHI}) and so, from \cite{DiBenedetto1983:art}[Theorem 1, Theorem 2], the sequence $\{\Phi_j\}_{j\in\NN}$ is uniformly bounded in $C^{1,\alpha}(\Omega')$ where $\Omega'$ is a any fixed compactly embedded subset of $\Omega$. Consequently, by the Ascoli-Arzel\`a theorem, we conclude that, up to passing to a suitable subsequence, $\nabla \Phi_j \to \nabla \underline{\Phi}$ uniformly in $\Omega'$ as $j \to +\infty$. Hence, for any $\phi \in \mathcal{C}_c^{\infty}(\Omega)$
\[
\int_{\Omega_j} |\nabla \Phi_j|^{p-2} \nabla\Phi_j \cdot \nabla \phi \,dx \to \int_{\Omega} |\nabla \underline{\Phi}|^{p-2} \nabla \underline{\Phi} \cdot \nabla \phi \,dx, \quad \text{ as } j \to +\infty,
\]
which shows that $\underline{\Phi} = \underline{\Phi}(x)$ is a weak solution to problem \eqref{eq:STATIONARYPROBLEM}, i.e., it satisfies \eqref{eq:EQUATIONFORPHIJ} for any $\phi \in \mathcal{C}_c^{\infty}(\Omega)$. We are left to show that $\underline{\Phi} = \underline{\Phi}(x)$ does not depend on $y \in \RR$. If we do so, the thesis follows by using the uniqueness of $\Phi = \Phi(z)$, inequality \eqref{eq:INEQUALITYOFLIMITSFUNCTIONPHI}, and the fact that the set $K \subset \Omega$ is bounded (note that here we strongly use the regularity w.r.t. $x \in \RR^{N+1}$ of the sequence $\Phi_j = \Phi_j(x)$). \normalcolor So, for any $a > 0$, let us consider the sequence $\Phi_j^a = \Phi_j^a(x)$ of the unique nonnegative weak solutions to the stationary problem \eqref{eq:STATIONARYPROBLEM}:
\[
\begin{cases}
-(p-2)\Delta_p \Phi_j^a = \Phi_j^a \quad &\text{ in } \Omega_j^a \\
\Phi_j^a = 0              \quad &\text{ in } \partial \Omega_j^a,
\end{cases}
\]
where, we have defined:
\[
\Omega_j^a = D \times(-j+a,j+a),
\]
for any integer $j \geq 1$, and the sequence $\Phi_{j+a} = \Phi_{j+a}(x)$ of the unique nonnegative weak solutions to the same problem posed in
\[
\Omega_{j+a} = D \times(-j-a,j+a),
\]
for any integer $j \geq 1$. Note that by uniqueness it immediately follows $\Phi_j^a(z,y) = \Phi_j(z,y-a)$ and, since $(-j+a,j+a) \subset (-j-a,j+a)$ we have $\Omega_j^a \subset \Omega_{j+a}$, so that
\[
\Phi_j^a \leq \Phi_{j+a} \quad \text{ in } \Omega_j^a,
\]
by comparison. Consequently, passing to the limit as $j \to +\infty$ in the above inequality, we obtain
\[
\underline{\Phi}(z,y-a) \leq \underline{\Phi}(z,y) \quad \text{ in } \Omega,
\]
for any $a > 0$. Finally, repeating this procedure for $a < 0$, we get that above relation is indeed an equality and we conclude the proof by the arbitrariness of $a > 0$. $\Box$

\paragraph{II. A special class of wave profiles.} We now prove a second preliminary result. As we have done above, we pass to the transformed solutions \eqref{eq:RENORMALIZATIONFORMULA} and, more precisely, we consider wave solutions to \eqref{eq:EQUATIONREACTIONTRANSFORMATION} which are independent of the variable $z \in D$, i.e. solutions in the form
\[
v(z,y,\tau) = \phi(\xi), \qquad \xi = y - c\tau,
\]
where $c > 0$ is the wave speed. We notice that this kind of travelling waves satisfy problem \eqref{eq:DIRICHLESTATIONARYTWS} replacing the homogeneous Dirichlet conditions with the Neumann ones: $\partial_{\nu} \phi = 0$ in $\partial\Omega$.

\noindent In this case, their analysis can be taken back to an ODEs problem in a nonstandard phase-plane. Indeed, plugging the above ansatz into \eqref{eq:REACTIONTRANSFORMATION}, we are led to the study of the degenerate second-order ODE:
\begin{equation}\label{eq:EQUATIONPROFILENEUMANN}
c\phi' + (|\phi'|^{p-2}\phi')' + \frac{\phi}{p-2} = 0, \qquad \xi \in \RR,
\end{equation}
with the convention $\phi' = d\phi/d\xi$.
\begin{lem}\label{Lemma:ExistenceTWNeumann}
For any $c > 0$, there exists a unique solution $\phi_c = \phi_c(\xi)$ to \eqref{eq:EQUATIONPROFILENEUMANN} such that
\[
\phi_c(\xi) = 0 \quad \text{ for all } \xi \geq \xi_0,
\]
for some $\xi_0 \in \RR$ (depending on $c > 0$). Uniqueness is understood up to horizontal shifts. Moreover, the function
\[
c \to M_c := \max_{\xi \in \RR} \phi_c(\xi)
\]
is well-defined in $(0,+\infty)$, monotone nondecreasing, with $M_c \to +\infty$ as $c \to +\infty$.
\end{lem}
\paragraph{Proof.} Assuming $\phi \geq 0$ and following \cite{DePablo-Vazquez1:art} and \cite{AA2017:art,AA-JLV:art}, we introduce the new variables
\[
X = \phi \qquad \text{ and } \qquad Z = -\left(\frac{p-1}{p-2}\,\phi^{\frac{p-2}{p-1}}\right)' = -X^{-\frac{1}{p-1}}X',
\]
which stand for the density and the derivative of the pressure profile (cf.  \cite{EstVaz:art}). We thus obtain the first-order ODE system
\begin{equation}\label{eq:SYSTEMNONSINGULARTWSXI}
-\frac{dX}{d\xi} = X^{\frac{1}{p-1}}Z, \quad\quad -(p-1)X^{\frac{p-2}{p-1}}|Z|^{p-2}\, \frac{dZ}{d\xi} = cZ - |Z|^p - \frac{X^{\frac{p-2}{p-1}}}{p-2},
\end{equation}
which, after re-parametrization $d\xi = -(p-1)X^{\frac{p-2}{p-1}}|Z|^{p-2}ds$, takes the nonsingular form
\[
\frac{dX}{ds} = (p-1)X|Z|^{p-2} Z, \quad\quad \frac{dZ}{ds} = cZ - |Z|^p - \frac{X^{\frac{p-2}{p-1}}}{p-2}.
\]
Of course, the two systems are equivalent outside their critical points $O(0,0)$ and $R_c(0,c^{1/(p-1)})$. The \emph{equation of the trajectories} is
\begin{equation}\label{eq:EQUATIONTRAJECTORIESTWS}
\frac{dZ}{dX} = \frac{(p-2)(cZ - |Z|^p) - X^{\frac{p-2}{p-1}}}{(p-1)(p-2)X|Z|^{p-2}Z} := H(X,Z;c),
\end{equation}
obtained by eliminating the time variable. As first observation, we note that $H(X,Z;c) > 0$ if $X > 0$ and $Z < 0$, for any $c \geq 0$. This fact, combined with the uniqueness of solutions at regular points, implies that any trajectory $Z = Z(X)$ defined for $X > 0$ satisfies $Z(X) \sim -\infty$ for $X \sim 0^+$. Consequently, from \eqref{eq:EQUATIONTRAJECTORIESTWS}, it must be
\begin{equation}\label{eq:ASBEHCHANGESIGNPOINTS}
(p-1)\frac{dZ}{dX} \sim -\frac{Z}{X} \quad \text{ for } X \sim 0^+ \qquad \text{ i.e. } \qquad Z(X) \sim - X^{-\frac{1}{p-1}} \quad \text{ for } X \sim 0^+.
\end{equation}
Now, moving to the subset of the phase plane where $X,Z > 0$, it is not hard to see that the \emph{null isoclines} are composed by a unique branch (belonging to the first quadrant), linking the points $O$, $R_c$, and recalling a horizontal parabola with vertex in
\begin{equation}\label{eq:EXPRESSIONVERTEXNULLISOCLINES}
(X_c,Z_c) = \left([(p-1)(p-2)]^{\frac{p-1}{p-2}}\left(\frac{c}{p}\right)^{\frac{p}{p-2}}, (p-1)\left(\frac{c}{p}\right)^{\frac{p}{p-1}} \right).
\end{equation}
The shape of the \emph{null-isoclines} and the study of the sign of \eqref{eq:EQUATIONTRAJECTORIESTWS} show that for any $c > 0$ there are two classes of special trajectories:

(i) There is exactly one trajectory $T_c = T_c(X)$ coming from the saddle-type point $R_c(0,c^{1/(p-1)})$ which is monotone decreasing when $Z > 0$ and intersect the axis $Z = 0$ in a point $X := M_c > 0$ (cfr. with \cite{AA2017:art,AA-JLV:art}). It satisfies the Darcy's law in its free boundary point $\xi_0 \in \RR$:
\[
\phi(\xi_0)^{-\frac{1}{p-1}}\phi'(\xi_0) = - c^{\frac{1}{p-1}} \qquad \text{ i.e. } \qquad \phi^{\frac{p-2}{p-1}}(\xi) \sim \frac{p-2}{p-1} c^{\frac{1}{p-1}} (\xi_0 - \xi), \quad \text{ for } \xi \sim \xi_0.
\]
This fact can be easily verified by noting that $T_c(X) \sim c^{\frac{1}{p-1}}$ for $X \sim 0^+$, so that the first equation in \eqref{eq:SYSTEMNONSINGULARTWSXI} satisfies $-X' \sim (c X)^{1/(p-1)}$ for $X \sim 0^+$.
This kind of wave profiles are known in literature as \emph{fast} or \emph{finite} orbits, and they satisfies $\phi(\xi) = 0$ for any $\xi \geq \xi_0$. Note that the value $M_c > 0$ is nothing more than the maximum value assumed by $\phi = \phi(\xi)$ in the all $\RR$ (this easily follows from the fact that $dZ/dX > 0$ for $X > 0$ and $Z < 0$) and, furthermore, the following bound holds (cfr. with \eqref{eq:EXPRESSIONVERTEXNULLISOCLINES}):
\begin{equation}\label{eq:BOUNDBELOWMAXPOINTFASTTW}
M_c \geq X_c := [(p-1)(p-2)]^{\frac{p-1}{p-2}}\left(\frac{c}{p}\right)^{\frac{p}{p-2}} \quad \text{ for all } c > 0.
\end{equation}
Consequently, $M_c \to + \infty$ as $c \to +\infty$. Finally, note that plugging \eqref{eq:ASBEHCHANGESIGNPOINTS} into the first equation in \eqref{eq:SYSTEMNONSINGULARTWSXI}, we obtain $\phi' \sim 1$ for $X \sim 0$, $Z \sim -\infty$, which proves that the profile $\phi = \phi(\xi)$ changes sign in finite time.

\begin{figure}[h!]
  \centering
  \includegraphics[scale=0.5]{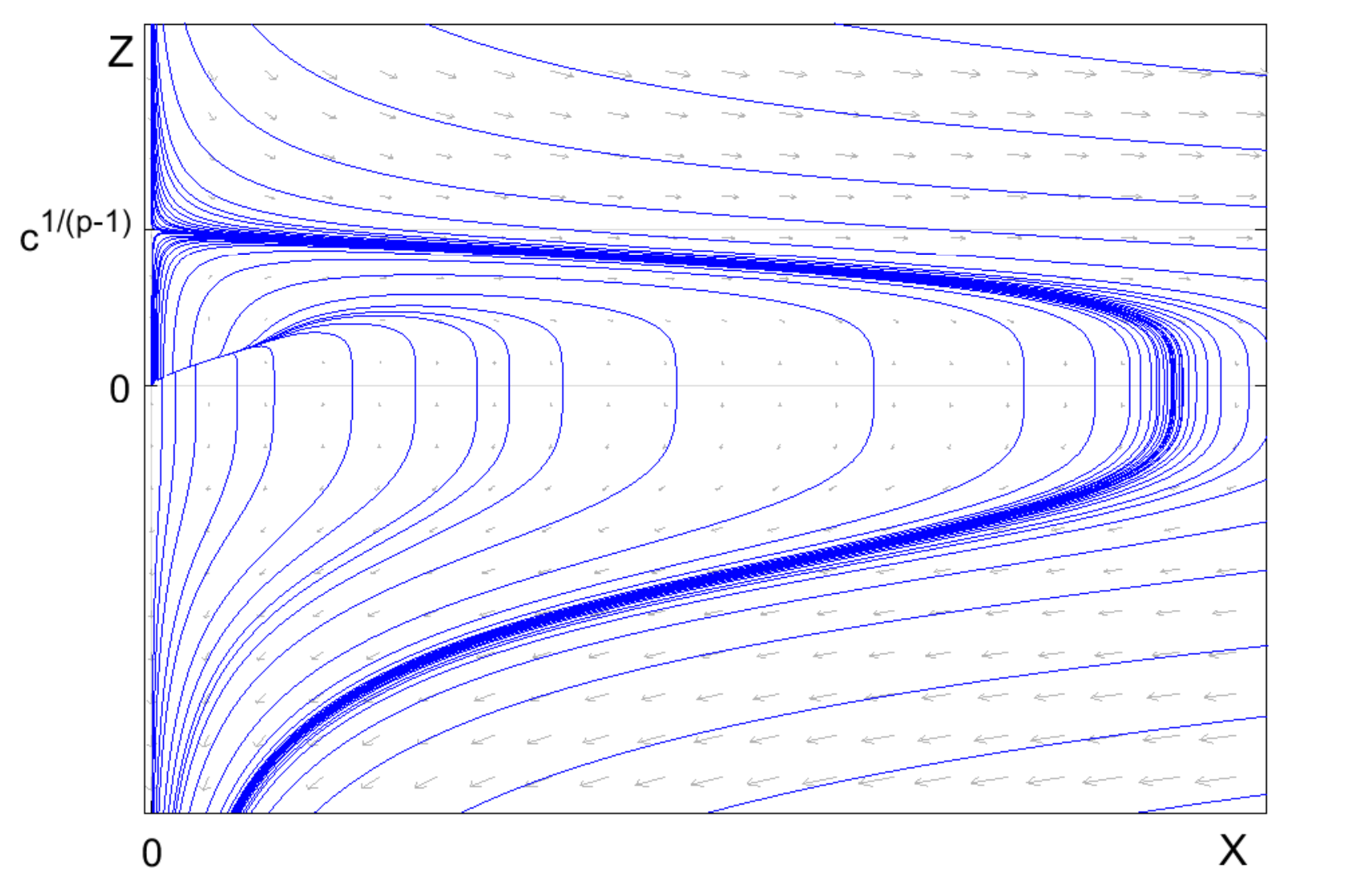}
  \caption{A qualitative representation of the trajectories in the $(X,Z)$-phase plane.}\label{Fig:NeumannWaves}
\end{figure}

(ii) There are infinitely many trajectories $Z = Z(X)$ from the (unstable) node type point $O(0,0)$ above the null-isocline branch, given by the curve $\widetilde{Z}(X) = 1/[c(p-2)] X^{(p-2)/(p-1)}$ for $X \sim 0^+$. Their asymptotic behaviour near the point $O(0,0)$ has not an immediate analytic expression, but, as the \emph{fast} one have a unique maximum point and change sign in finite time. We do not insist on these orbits since they will not be important in the following.

\emph{Monotonicity of $T_c = T_c(X)$ w.r.t. $c > 0$.} We are left to show that the function $c \to M_c$ is monotone nondecreasing. We will prove that if $T_c = T_c(X)$ denotes the fast orbit from $R_c$, then
\[
\text{if } \quad 0 < c_1 < c_2 \quad \text{ then } \quad T_{c_1}(X) < T_{c_2}(X) \quad \text{ for any } 0 \leq X < M_{c_1},
\]
from which our statement follows. First of all, we note that
\begin{equation}\label{eq:MONOTONICITYRESPECTCEQTRAJECTORIESTYPEC}
\frac{\partial H}{\partial c}(X,Z;c) = \frac{1}{(p-1)X|Z|^{p-2}} > 0, \quad \text{for all }  X,Z > 0, \; c > 0,
\end{equation}
which implies $H(X,Z;c_1) < H(X,Z;c_2)$ for any $c_1 < c_2$. Now, assume by contradiction that $T_{c_1}$ and $T_{c_2}$ touch in a point $(X_0,T_{c_1}(X_0) = T_{c_2}(X_0))$, with $0 < X_0 < M_{c_1}$. Up to take a smaller contact point, we can assume $T_{c_2} > T_{c_1}$ for any $0 \leq X < X_0$. Consequently, for $h > 0$ small enough, it must be
\[
\frac{T_{c_2}(X_0) - T_{c_2}(X_0 - h)}{h} \leq \frac{T_{c_1}(X_0) - T_{c_1}(X_0 - h)}{h}
\]
and taking the limit as $h \to 0$, we get $dT_{c_2}(X_0)/dX \leq dT_{c_1}(X_0)/dX$, which is in contradiction with \eqref{eq:MONOTONICITYRESPECTCEQTRAJECTORIESTYPEC}. This ends the proof of the lemma. $\Box$

\paragraph{Remark.} We notice that we could have looked for wave solutions with the form
\[
v(z,y,\tau) = \phi(\xi), \qquad \xi = y + c\tau,
\]
where $c > 0$ and moving towards the left. The equation of the profile is
\[
-c\phi' + (|\phi'|^{p-2}\phi')' + \frac{\phi}{p-2} = 0, \qquad \xi \in \RR,
\]
and differs from \eqref{eq:EQUATIONPROFILENEUMANN} just for the minus sign in front of $c\phi'$. The change of variables
\[
X = \phi \qquad \text{ and } \qquad Z = \left(\frac{p-1}{p-2}\,\phi^{\frac{p-2}{p-1}}\right)',
\]
and the re-parametrization $d\xi = (p-1)X^{\frac{p-2}{p-1}}|Z|^{p-2}ds$ leads us to the same equation of the trajectories \eqref{eq:EQUATIONTRAJECTORIESTWS}. Consequently, the statement of Lemma \ref{Lemma:ExistenceTWNeumann} holds for this different kind of waves, but this time the profiles satisfy $\phi_c(\xi) = 0$ for all $\xi \leq \xi_0$ and moves towards the left.

%
%
%
%
%
%
%
%
%
%
%
\section{Existence, uniqueness and regularity of a special wave solution}\label{Section:ExistenceUniquenessTW}
In our long-time behaviour analysis of solutions to the transformed problem \eqref{eq:REACTIONTRANSFORMATION}, we will employ barriers built-up by using suitable wave solutions. More precisely, we will need solutions to problem \eqref{eq:REACTIONTRANSFORMATION} in the form
\begin{equation}\label{eq:GENERALTWFRONT}
v(z,y,\tau) = \varphi(z,y-c\tau) \qquad (z,y) \in \Omega, \; \tau > 0,
\end{equation}
which travel along the tube $\Omega$ with constant speed of propagation $c > 0$. This means that the wave \emph{profile} must satisfy the stationary Dirichlet problem
\begin{equation}\label{eq:DIRICHLESTATIONARYTWS}
\begin{cases}
\Delta_p\varphi + c \, \partial_{\xi} \varphi + \frac{\varphi}{p-2} = 0 \quad &\text{ in } \Omega \\
\varphi = 0                                                              \quad &\text{ in } \partial\Omega,
\end{cases}
\end{equation}
where $\xi = y -c\tau$ is the \emph{moving coordinate}. One of the main goals of this paper is to prove that such travelling wave profiles exist for some speed $c > 0$ and, at same time, connect the special solution $\Phi = \Phi(z)$ to \eqref{eq:STATIONARYPROBLEM} (at $\xi = -\infty$) to the level $v = 0$ (at $\xi = +\infty$), namely
\begin{equation}\label{eq:LIMITCONDTWPROFILES}
\lim_{\xi \to -\infty} \varphi(z,\xi) = \Phi(z), \qquad \lim_{\xi \to +\infty} \varphi(z,\xi) = 0, \quad z \in D,
\end{equation}
which are nothing more than \eqref{eq:LIMITCONDTWPROFILES0}. We first prove the existence of such solutions together with some their qualitative properties and then we discuss the problem of uniqueness, which turns out to be less difficult problem.

\paragraph{Existence and regularity.} This first part is devoted to prove the existence of a wave solution to problem \eqref{eq:DIRICHLESTATIONARYTWS}-\eqref{eq:LIMITCONDTWPROFILES}. Furthermore, we show some regularity properties of its profile and the finiteness of its free boundary. In more precise words, we prove the following lemma.
\begin{lem}\label{Theorem:WAVECONSTRUCTION}
There exists a speed $c_{\ast} > 0$ for which problem \eqref{eq:DIRICHLESTATIONARYTWS}-\eqref{eq:LIMITCONDTWPROFILES} has a weak solution $\varphi = \varphi(z,\xi)$ which is continuous and monotone non-increasing w.r.t. the longitudinal variable $\xi$.
\\
Furthermore, $\varphi = \varphi(z,\xi)$ has a free boundary which is a bounded, connected and locally smooth subset of $\Omega$ with Hausdorff codimension one (in $\RR^{N+1}$), which can be locally parametrized by a locally smooth function.
\end{lem}
\paragraph{Proof.} We notice that, once the existence of a TW with free boundary is proven, its properties and the ones of its free boundary will follow by the well-known regularity theory and parabolic $p$-Laplacian equations. Now, we fix $c > 0$ and we construct our special solution through the following approximation procedure.

\emph{Step 1: Approximating sequence.} To approximate solutions to problem \eqref{eq:DIRICHLESTATIONARYTWS}-\eqref{eq:LIMITCONDTWPROFILES}, we write $v = v(z,y,\tau)$ w.r.t. the moving coordinate system, i.e., we consider the change of variable
\[
v(z,y,\tau) = w(z,\xi,\tau), \quad \xi = y - c\tau,
\]
so that the equation in \eqref{eq:REACTIONTRANSFORMATION} is transformed into
\[
\partial_{\tau} w = \Delta_pw + c\partial_{\xi}w + \frac{w}{p-2} \quad \text{ in } \Omega\times(0,\infty).
\]
More precisely, we consider a sequence of weak solutions $w_j = w_j(z,y,\tau)$ to the problems
\begin{equation}\label{eq:APPROXTRANSFPROB}
\begin{cases}
\partial_{\tau} w_j = \Delta_pw_j + c\partial_{\xi}w_j + \frac{w_j}{p-2} \quad &\text{ in } \Omega_j\times(0,\infty) \\
w_j = b_j                 \quad &\text{ in } \partial\Omega_j\times(0,\infty) \\
w_j(\cdot,\cdot,0) = 0       \quad &\text{ in } \Omega_j,
\end{cases}
\end{equation}
where $\Omega_j = D\times(-j,j)$,
\[
b_j :=
\begin{cases}
\Phi \quad &\text{ in } D\times\{y = -j\} \\
0    \quad &\text{ in } \partial\Omega_j \setminus (D\times\{y = -j\}).
\end{cases}
\]
and $\Phi = \Phi(z)$ is the unique nonnegative weak solution to the stationary problem \eqref{eq:STATIONARYPROBLEM}. Solutions to problem \eqref{eq:DIRICHLESTATIONARYTWS}-\eqref{eq:LIMITCONDTWPROFILES} will be obtained as double limit $\tau \to +\infty$ and $j \to +\infty$ of the approximating sequence $w_j = w_j(z,y,\tau)$.

\noindent As a first observation, note that the combined use of the maximum principle and the positivity of the stationary super-solution $\Phi = \Phi(z)$, immediately implies that
\[
0 < w_j(z,\xi,\tau) \leq \Phi(z), \quad \text{ in } \Omega_j\times(0,\infty),
\]
for any $j \in \NN$. We notice that from the above inequality we deduce that solutions to \eqref{eq:APPROXTRANSFPROB} exist for any time $\tau > 0$. Moreover, they are monotone non-increasing w.r.t. to the longitudinal variable and non-decreasing w.r.t. time variations:
\[
\partial_{\xi} w_j \leq 0 \quad \text{ in } \Omega_j\times(0,\infty), \qquad \partial_{\tau} w_j \geq 0 \quad \text{ in } \Omega_j\times(0,\infty).
\]
Both inequalities can be proved through comparison techniques. To show the first one, for any $h > 0$, we introduce the translated $w_j^h(z,\xi,\tau) := w_j(z,\xi + h,\tau)$, and we note that it satisfies the problem
\[
\begin{cases}
\partial_{\tau} w_j^h = \Delta_pw_j^h + c\partial_{\xi}w_j^h + \frac{w_j^h}{p-2} \quad &\text{ in } \Omega_j^h\times(0,\infty) \\
w_j^h \leq w_j                 \quad &\text{ in } \partial\Omega_j^h\times(0,\infty) \\
w_j^h(\cdot,\cdot,0) = 0       \quad &\text{ in } \Omega_j^h,
\end{cases}
\]
where $\Omega_j^h := D\times(-j,j-h)$. The inequality on the boundary holds since $w_j^h(z,j-h,\tau) = 0$ and $w_j^h(z,\xi,\tau) \leq \Phi(z)$ for any $z \in D$. Consequently, we get $w_j^h \leq w_j$ in $\Omega_j^h\times(0,\infty)$ which gives the desired inequality thanks to the arbitrariness of $h > 0$. For what concerns the second one, it is enough to note that, since $w_j \geq 0$, then translation $w_j^s(z,\xi,\tau) := w_j(z,\xi,\tau + s)$ ($s > 0$) is a super-solution to problem \eqref{eq:APPROXTRANSFPROB} and so, using the arbitrariness of $s > 0$ we conclude the proof of the above claim.

\emph{Step 2: Limit as $\tau \to +\infty$.} The time-monotonicity and uniform upper-bound $0 \leq w_j = w_j(z,\xi,\tau) \leq \Phi(z)$ imply that for any $j \in \NN$:
\[
w_j(\cdot,\cdot,\tau) \to \varphi_j(\cdot,\cdot) \quad \text{ in } \Omega_j,
\]
as $\tau \to +\infty$ point-wise, for some nontrivial and nonnegative limit function $\varphi_j = \varphi_j(z,\xi)$, satisfying also
\[
0 \leq \varphi_j(z,\xi) \leq \Phi(z), \qquad \partial_{\xi}\varphi_j \leq 0 \quad \text{ in } \Omega_j,
\]
for any $j \in \NN$.  Since the above convergence is monotone non-decreasing in $\tau > 0$, it must be $\partial_{\tau}w_j(z,\xi,\tau) \to 0$ point-wise in $\Omega_j$ as $\tau \to +\infty$ (this easily follows by showing that $\partial_{\tau}w_j(z,\xi,\tau)$ is bounded for any $(z,\xi) \in \Omega_j$ as $\tau \to +\infty$). Consequently, using again the regularity theory for $p$-Laplacian equations  \cite{DiBenedetto1983:art}[Theorem 1, Theorem 2], we get that the limit $\varphi_j =\varphi_j(z,\xi)$ is a (bounded and $C^{1,\alpha}(\Omega_j)$) weak solution to
\begin{equation}\label{eq:APPROXTRANSFPROBSTATIONARY}
\begin{cases}
\Delta_p\varphi_j + c \partial_{\xi} \varphi_j + \frac{\varphi_j}{p-2} = 0 \quad &\text{ in } \Omega_j \\
\varphi_j = b_j                                                      \quad &\text{ in } \partial\Omega_j,
\end{cases}
\end{equation}
it is strictly positive, and $\varphi_j < \Phi$ in $\Omega_j$ thanks to the strong maximum principle. \normalcolor Actually, we can show it is the \emph{unique} weak solution to \eqref{eq:APPROXTRANSFPROBSTATIONARY} through the sliding method as follows. Assume there is another solution $\psi_j = \psi_j(z,\xi)$ and take $h > 0$ large enough such that
\[
\psi_j(z,\xi) \leq \varphi_j(z,\xi - h) \quad \text{ in } D\times(-j+h,j).
\]
The validity of the above inequality comes from the regularity of our solutions and their values at the boundary of $D\times(-j+h,j)$ (once $h > 0$ is suitably chosen). Now, sliding back to the left, we pick $\underline{h} > 0$ defined as
\[
\underline{h} := \min\left\{h \geq 0: \psi_j(z,\xi) \leq \varphi_j(z,\xi - h), \quad (z,\xi) \in D\times(-j+h,j) \right\} > 0.
\]
Now, since both $\varphi_j$ and $\psi_j$ are strictly positive in $\Omega_j$ we obtain that $\psi_j(z,\xi) \leq \varphi_j(z,\xi-\underline{h})$ and touch at a point belonging to $D\times(-j+\underline{h},j)$, which is in contradiction with the strong maximum principle.

\emph{Step 3: Monotonicity and continuity w.r.t. $c > 0$.} Since $\partial_{\xi} \varphi \leq 0$, an easy comparison procedure shows that if $\varphi_{j,c} = \varphi_{j,c}(z,\xi)$ denotes the unique solution to \eqref{eq:APPROXTRANSFPROBSTATIONARY}, then the function $c \to \varphi_{j,c}$ is monotone non-increasing, in the sense that if $c_1 \leq c_2$, then $\varphi_{j,c_1} \geq \varphi_{j,c_2}$ in $\Omega_j$.  Hence, for any $j\in \NN$ and $c_0 > 0$, we know that the point-wise limit
\[
\widetilde{\varphi}_j := \lim_{c \to c_0}\varphi_{j,c} \quad \text{ in } \Omega_j
\]
exists (it follows by the monotonicity above and the bound on $\varphi_{j,c}$). Moreover, using the regularity of $\varphi_{j,c}$ and its limit, we get that $\widetilde{\varphi}_j = \widetilde{\varphi}_j(z,\xi)$ is a weak-solution to \eqref{eq:APPROXTRANSFPROBSTATIONARY} with $c = c_0$ and so, by uniqueness, it must be $\widetilde{\varphi}_j = \varphi_{j,c_0}$, i.e. the function $c \to \varphi_{j,c}$ is continuous (for any $j \in \NN$). \normalcolor

\emph{Step 4: Passage to the limit as $j \to +\infty$.} The next step is to pass to the limit as $j \to +\infty$ in our approximating sequence. The main problem here is to show that $\varphi_j \not \to 0$ or $\varphi_j \not \to \Phi$ uniformly as $j \to +\infty$. These facts will be proved by properly choosing $c = c_j$ for large $j \in \NN$.

\noindent CLAIM: For any $j \in \NN$ large, there is a speed $c_{j\ast} > 0$ such that
\begin{equation}\label{eq:NONTRIVIALITYOFTHELIMITTHM}
\varphi_{j,c_{j\ast}}(0,0) = \frac{1}{2} \Phi(0).
\end{equation}
Furthermore, there are two constants $0 < \underline{c} < \overline{c}$ such that $\underline{c} < c_{j\ast} < \overline{c}$, for any $j \in \NN$ large enough.

To prove the claim we proceed in some steps as follows.

(i) First of all, we note that if $c = 0$, we have $\varphi_{j,0}(z,\xi) \to \Phi(z)$ uniformly on compact sets of $\Omega$. This follows by comparing $\varphi_{j,0} = \varphi_{j,0}(z,\xi)$ with the sequence of solutions $\Phi_j = \Phi_j(z,\xi)$ to problem \eqref{eq:APPROXSEPVARSOLLEMMACONVONCOMPACT} (cfr. with the proof of Lemma \ref{LEMMACONVERGENCEINBOUDEDDOMAINS}) and recalling that $\varphi_{j,0}(z,\xi) \leq \Phi(z)$ for any $j \in \NN$ by construction.
\\
Therefore, there is $j_1 \in \NN$ large enough such that $\varphi_{j,0}(0,0) > \Phi(0)/2$ for any $j \geq j_1$, and so, by the continuity of the map $c \to \varphi_{j,c}$, we get
\[
\varphi_{j,c}(0,0) > \frac{1}{2}\Phi(0) \quad \text{ for any } 0 < c < \underline{c}_0,
\]
for some suitable $\underline{c}_0$ (depending on $j \in \NN$). Secondly, we show that there exist $j_2 \in \NN$ and $\overline{c}_0 > 0$, such that
\begin{equation}\label{eq:LIMITVARPHILARGEC}
\varphi_{j,c}(0,0) < \frac{1}{2}\Phi(0) \quad \text{ for any } c > \overline{c}_0, \; j \geq j_2.
\end{equation}
This is crucial. Note indeed that once it is proved, \eqref{eq:NONTRIVIALITYOFTHELIMITTHM} follows by the monotonicity and the continuity of the function $c \to \varphi_{j,c}$, once we take $j \geq \max\{j_1,j_2\}$. The proof of \eqref{eq:LIMITVARPHILARGEC} is not trivial and we devote to it a separate step.

(ii) The idea is to compare solutions $w_{j,c} = w_{j,c}(z,\xi,\tau)$ to \eqref{eq:APPROXTRANSFPROB} with the dynamic version of \emph{fast} TWs $\phi_c = \phi_c(\xi)$ found in Lemma \ref{Lemma:ExistenceTWNeumann}. We recall that each $\phi_c$ has a maximum $M_c \to +\infty$ as $c \to + \infty$ and bounded support to its right.

Now, let us fix $\overline{c}_0 > 0$ such that $\phi_{\overline{c}_0}(0) = M_{\overline{c}_0} \geq \Phi(z)$ for any $z \in D$, and take $c > \overline{c}_0$ (this is possible since $\Phi = \Phi(z)$ is bounded in $D$). Furthermore, if $\phi_{\overline{c}_0} = \phi_{\overline{c}_0}(s)$ denotes the \emph{fast} TW with speed of propagation $\overline{c}_0 > 0$ (and $s := y - \overline{c}_0\tau$), we set $c_0 := c-\overline{c}_0 > 0$ and we define
\[
\psi^l(z,\xi,\tau) := \phi_{\overline{c}_0}(\xi + c_0\tau-l), \quad \text{ where } \xi = y - c\tau,
\]
and $l \in [-j,j])$. Note that (up to a shift) we can assume $\phi_{\overline{c}_0}(0) = M_{\overline{c}_0}$ and, since $c_0 > 0$, the profile of $\psi^l = \psi^l(z,\xi,\tau)$ travels towards the left. We want to employ it as super-solution to $w_{j,c} = w_{j,c}(z,\xi,\tau)$. To do so, we introduce the following comparison domain:
\[
Q_{j,l} := \left\{(z,\xi,\tau) \in \Omega_j\times(0,\infty) : \xi + c_0\tau \geq l, \; 0 \leq \tau < \tau_{j,l}\right\}, \quad \text{ where } \; \tau_{j,l}: = \frac{l+j}{c_0}.
\]
The choice of the sub-domain $Q_{j,l}$ allows us to compare $w_{j,c}$ and $\psi_l$ on the boundary. Indeed, we have
\[
w_{j,c}(z,\xi,\tau) \leq \Phi(z) \leq \phi_{\overline{c}_0}(0) \equiv \psi^l(z,\xi,\tau) \quad \text{ in } \{(z,\xi,\tau) \in \overline{Q}_{j,l}: \xi + c_0\tau = l, \;  0 \leq \tau \leq \tau_{j,l}\},
\]
up to take $\overline{c}_0 > 0$ larger (this is possible since $M_c \to +\infty$ as $c \to +\infty$). In the other part of the boundary, the comparison is trivial since $w_{j,c} \equiv 0$, while $\psi^l \geq 0$. Similar for the comparison at the initial time $\tau = 0$. We thus have $w_{j,c}(z,\xi,\tau) \leq \psi^l(z,\xi,\tau)$ in $Q_{j,l}$ and, in particular,
\[
w_{j,c}(z,\xi,\tau_{j,l}) \leq \psi^l(z,\xi,\tau_{j,l}) = \phi_{\overline{c}_0}(\xi + j) \quad \text{ in } \Omega_j.
\]
Consequently, since $l \in [-j,j]$ is arbitrary we can take $\tau_{j,l} \in [0,2j/c_0]$, and so
\[
w_{j,c}(z,\xi,\tau) \leq \psi^l(z,\xi,\tau) = \phi_{\overline{c}_0}(\xi + j) \quad \text{ in } \Omega_j\times[0,2j/c_0],
\]
which implies
\[
w_{j,c}(0,0,\tau) \leq \phi_{\overline{c}_0}(j) \quad \text{ for any } 0 \leq \tau \leq 2j/c_0, \; j \in \NN.
\]
Consequently, there is $j_2 > 0$ large enough (depending on $\overline{c}_0 > 0$) such that $\phi_{\overline{c}_0}(j) = 0$ for any $\RR \ni j \geq j_2$ and so
\[
w_{j,c}(0,0,\tau) = 0 \quad \text{ for any } 0 \leq \tau \leq 2j/c_0, \; j \geq j_2.
\]
This complete the proof of \eqref{eq:LIMITVARPHILARGEC} since $w_{j,c}(0,0,\tau) \to \varphi_{j,c}(0,0)$ as $\tau \to +\infty$ (taking eventually $j_2 \in \NN$ larger).

\emph{Step 5: Construction of the TW.} Now, let $j \in \NN$ be large and $0 < \underline{c} < c_j := c_{\ast j} < \overline{c}$ such that \eqref{eq:NONTRIVIALITYOFTHELIMITTHM} holds. Up to passing to a subsequence, we can assume
\[
c_j \to c_{\ast} \quad \text{ as } j \to +\infty, \quad \text{ with } \underline{c} \leq c_{\ast} \leq \overline{c}.
\]
 Moreover, up to passing to another subsequence, we have (using the usual regularity estimates)
\[
\varphi_{j,c_j} \to \varphi \quad \text{ in } \mathcal{C}^{1,\alpha}(\Omega'),
\]
for any compact set $\Omega' \subset \Omega$ and for some continuous function $\varphi = \varphi(z,\xi)$ satisfying problem \eqref{eq:DIRICHLESTATIONARYTWS} with $c = c_{\ast}$ in the weak sense, together with the properties
\[
0 \leq \varphi(z,\xi) \leq \Phi(z) \quad \text{ and } \quad \partial_{\xi}\varphi(z,\xi) \leq 0 \quad \text{ in } \Omega.
\]
Note that thanks to \eqref{eq:NONTRIVIALITYOFTHELIMITTHM} we can assume both $\varphi \not= 0$ and $\varphi \not= \Phi$. Finally, note that $\varphi(z,\xi) \to \Phi(z)$ as $\xi \to -\infty$ for any $z \in D$. Indeed, the point-wise limit
\[
\underline{\Phi}(z) := \lim_{\xi \to -\infty} \varphi(z,\xi),
\]
exists by monotonicity (w.r.t. $\xi \in \RR$) with $0 < \underline{\Phi} \leq \Phi$ and satisfies problem \eqref{eq:STATIONARYPROBLEM}. So, it must be $\underline{\Phi} = \Phi$ by uniqueness (of $\Phi = \Phi(z)$) and the proof of the first limit of \eqref{eq:LIMITCONDTWPROFILES} is completed. At the same time, it must be
\[
\lim_{\xi \to +\infty} \varphi(z,\xi) = 0, \quad \text{ for all } z \in D.
\]
If not, the limit function turns out to be $\Phi = \Phi(z)$ by the argument above. However, this is in contradiction with the fact that $\varphi \not= \Phi$ and $\varphi(\cdot,\xi) \to \Phi(\cdot)$ as $\xi \to -\infty$.
\normalcolor

\emph{Step 6: Bounds for the TW's free boundary.} If $\varphi = \varphi(z,\xi)$ denotes the TW just constructed, we define its free boundary by
\begin{equation}\label{eq:TWFREEBOUNDARY}
S := \{(z,S(z)): z \in D \} \quad \text{ where } \quad S(z) := \inf\{\xi \in \RR: \varphi(z,\xi) = 0 \}, \quad z \in D.
\end{equation}
Note that the above definition is well-posed thanks to the monotonicity w.r.t. $\xi \in \RR$ and, a priori, the set $S$ could be unbounded. The main goal of this fifth step is to show that this situation cannot happen. More precisely, we show that the wave's support is bounded to the right while the free boundary function $\xi = S(z)$ is bounded both below (left) and above (right). Recall that from the step above, we know that $\varphi$ has speed of propagation $c_{\ast} > 0$ and $\varphi(z,\xi) \to 0$ as $\xi \to +\infty$, uniformly in $z \in D$ (this easily follows by construction). Moreover, up to a shift (w.r.t to the variable $\xi$), we can assume $\varphi(0,0) = \Phi(0)/2$.

\noindent The idea is to proceed similar to \emph{Step 4} comparing the approximating sequence $w_{j,c_j} = w_{j,c_j}(z,\xi,\tau)$ with the \emph{fast} TWs constructed in the ODEs analysis. We have $c_j \to c_{\ast}$ as $j \to +\infty$ and we define
\[
\psi^l(z,\xi,\tau) := \phi_{c_{\ast}}(\xi - l),
\]
where $\phi_{c_{\ast}}$ is the \emph{fast} TW with speed $c_{\ast} > 0$, and $l > 0$ is chosen so that
\[
\max_{z \in D} \varphi(z,l) \leq \phi_{c_j}(0),
\]
for any $j \in \NN$ large enough (note that this is possible since $c_j \to c_{\ast}$ and thanks the continuity of the map $c \to \phi_c$). Now, we consider a new comparison domain
\[
Q_l := \{(z,\xi,\tau) \in \Omega\times(0,\infty): \xi > l \},
\]
and we note that by the monotonicity properties proved in the above steps, we have
\[
w_{j,c_j}(z,l,\tau) \leq \varphi_{j,c_j}(z,l) \leq \varphi(z,l) \leq \phi_{c_{\ast}}(0), \quad \text{ for any } z \in D, \; \tau > 0,
\]
and $j \in \NN$ large enough, where we recall that $c_j := c_{\ast j}$. Furthermore, since as always $w_{j,c_j} = 0$ in $\partial\Omega\times(0,\infty)$, it follows
\[
w_{j,c_j}(z,\xi,\tau) \leq \phi_{c_{\ast}}(\xi - l) \quad \text{ in } Q_l,
\]
and so, passing to the limit as $\tau \to +\infty$ and then as $j \to +\infty$, we easily get
\[
\varphi(z,\xi) \leq \phi_{c_{\ast}}(\xi - l) \quad \text{ in } Q_l,
\]
and since $\phi_{c_{\ast}}(\cdot)$ has bounded support to the right, $\varphi(\cdot)$ as bounded support to the right, too.

\noindent We are left to show that the free boundary $\xi = S(z)$ is bounded below (left). Assume by contradiction that there is $z_0 \in D$ such that $S(z_0) = -\infty$. Consequently, from the fact that $\partial_{\xi}\varphi \leq 0$, it must be $\varphi(z_0,\xi) = 0$ for any $\xi \in \RR$. However, using the fact that $\varphi_j(z_0,\xi) \to \varphi(z_0,\xi)$ for any $\xi \in \RR$ and $\varphi_j(z_0,-j) = \Phi(z_0) > 0$, we obtain the desired contradiction taking $j \in \NN$ large enough.

\noindent Finally, the properties of the free boundary follow by the standard regularity theory of solutions to equations with $p$-Laplacian diffusion. $\Box$

\paragraph{Important remarks.} We complete this subsection with some crucial remarks. First of all, in what follows, we will need to employ a strengthened version of the limit $\Phi(z) = \lim_{\xi \to -\infty}\varphi(z,\xi)$ proved in \emph{Step 5} of the above proof. Indeed, we have that $\varphi(\cdot,\xi)$ converges in \emph{relative error} to $\Phi(\cdot)$, i.e.
\begin{equation}\label{eq:CONVERGENCEINRELATIVERROR}
\lim_{\xi \to -\infty} \frac{\varphi(\cdot,\xi)}{\Phi(\cdot)} = 1 \quad \text{ uniformly in } D.
\end{equation}
Since both $\Phi = \Phi(z)$ and $\varphi = \varphi(z,\xi)$ are positive for any $z \in D$ and $\xi \sim -\infty$, the above limit is easily obtained if we restrict the convergence to any compact subset of $D$. To verify the validity of \eqref{eq:CONVERGENCEINRELATIVERROR} near the boundary $\partial D$, we use a barrier argument based on the validity of the Hopf lemma for elliptic $p$-Laplacian type equations. The main fact (cfr. for instance with \cite{ManfrediVespri1994:art,S-V1:art}) is that both $\Phi = \Phi(z)$ and $\varphi = \varphi(z,\xi)$ behave like
\[
\Phi(z) \asymp \dist(z,\partial D), \qquad \varphi(z,\xi) \asymp \dist(z,\partial D), \quad \text{ for } z \sim \partial D\,,
\]
as a consequence of Hopf principle. \normalcolor Consequently, $\Phi = \Phi(z)$ and $\varphi = \varphi(z,\xi)$ are comparable near the boundary of $\Omega$ and the bound $\varphi(z,\xi) \geq (1 - \varepsilon)\Phi(z)$ holds for each $z \in D$ and $\xi \ll 0$ negative enough.

Finally, we point out that the above proof gives the existence of another TW $\psi = \psi(z,\xi)$ satisfying problem \eqref{eq:DIRICHLESTATIONARYTWS} with symmetric conditions at the ends of the tube
\begin{equation}\label{eq:LIMITCONDTWPROFILESSYMMETRIC}
\lim_{\xi \to -\infty} \psi(z,\xi) = 0, \qquad \lim_{\xi \to +\infty} \psi(z,\xi) = \Phi(z), \quad z \in D.
\end{equation}
To see this it is enough to repeat the proof above by changing longitudinal variable $\xi \to -\xi$. Of course, this time, $\psi = \psi(z,\xi)$ is monotone non-decreasing w.r.t. to $\xi$.
%
%
%
%
%
%
%
%
%
%
%
%
\subsection{Uniqueness}
\begin{lem}\label{Lemma:MONOTONICITYOFCASTWRTDOMAINS}
Let $D_1 \subset D_2$ be two domains and $\varphi_1 = \varphi_1(z,\xi)$ and $\varphi_2(z,\xi)$ be two proper TW solutions with speeds $c_{1\ast}$ and $c_{2\ast}$. If $\varphi_1$ is a finite TW, then $c_{1\ast} \leq c_{2\ast}$.
\end{lem}
\paragraph{Remark.} A TW profile $\varphi = \varphi(z,\xi)$ is said to be \emph{proper}, if it is a continuous weak solution to \eqref{eq:DIRICHLESTATIONARYTWS}-\eqref{eq:LIMITCONDTWPROFILES}, it is monotone non-increasing w.r.t. the longitudinal variable $\xi$, and satisfies $0 \leq \varphi < \Phi$ in $\Omega$.
\paragraph{Proof.} Let us define $\Omega_i := D_i\times(0,\infty)$ and $
v_i(z,y,\tau) := \varphi_i(z,y - c_{i\ast}\tau)$ in $\Omega_i\times(0,\infty)$, for $i = 1,2$. Coming back to the real time variable and so, to solutions to \eqref{eq:PLETUBULARDOMAIN}, we consider
\begin{equation}\label{eq:TWSORIGINALUI}
u_i(z,y,t) := t^{-\frac{1}{p-2}} \varphi_i(z,y - c_{i\ast}\ln t),
\end{equation}
for $i = 1,2$, and we introduce the functions $u_{1h}(z,y,t) = u_1(z,y + h,t)$ and $u_2^{\theta}(z,y,t) = u_2(z,y,t-\theta)$ for any $0 < \theta < 1$ fixed and suitable $h > 0$.

\noindent We want to compare $u_{1h}$ and $u_2^{\theta}$ in $\Omega_1 \times (1,\infty)$. At the initial time $\tau = 0$, i.e. $t = 1$, we have
\begin{equation}\label{eq:INEQUALITYFORPERTURBATIONSUI}
(1 - \theta)^{-\frac{1}{p-2}} \varphi_2(z,y - c_{2\ast}\ln (1 - \theta)) \geq \varphi_1(z, y + h),
\end{equation}
for some $h > 0$ large enough. Indeed, we recall that $\varphi_2(z,-\infty) = \Phi_2(z) > \Phi_1(z) \geq \varphi_1(z,y + h)$ for any $(z,y) \in \Omega_1$ (cfr. with the end of the statement of Theorem \ref{ASYMPTOTICSINBOUNDEDDOMAINS} and the proof of Theorem \ref{Theorem:WAVECONSTRUCTION}). Consequently, taking $h > 0$ large enough and using the monotonicity of $\varphi_i = \varphi_i(z,\xi)$ w.r.t. $\xi \in \RR $, together with the fact that $\varphi_1 = 0$ for $\xi$ large enough,  we see that the above inequality is satisfied for any $0 < \theta < 1$. The comparison at the boundary is trivial. Hence by comparison we deduce $u_2^{\theta}(z,y,t) \geq u_{1h}(z,y,t)$ in $\Omega_1\times(1,\infty)$, which means
\[
(t - \theta)^{-\frac{1}{p-2}} \varphi_2(z,y - c_{2\ast}\ln (t - \theta)) \geq t^{-\frac{1}{p-2}} \varphi_1(z,y + h - c_{1\ast}\ln t)
\]
and, taking the limit as $t \to +\infty$ in the above inequality, we thus deduce that $c_{1\ast}$ cannot exceed $c_{2\ast}$. However, note that this procedure does not imply $c_{1\ast} < c_{2\ast}$, which, by the way, is not true (even if $D_1 \subset D_2$ with strict inequality), due to the fact that TWs are given up to a shift along the $\xi$ direction. $\Box$
\begin{lem}\label{Lemma:APPROXIMATINGCASTFROMINSIDE}
Let $D \subset \RR^N$ be a bounded domain and $D_1 \subset D_2 \subset \ldots \subset D$ a sequence of subsets such that
\[
D = \bigcup_{j \in \NN} D_j.
\]
Then there exists a sequence $\varphi_j = \varphi_j(z,\xi)$ of solutions to problem \eqref{eq:DIRICHLESTATIONARYTWS}-\eqref{eq:LIMITCONDTWPROFILES} posed in $\Omega_j = D_j\times\RR$ which converge to a solution $\varphi = \varphi(z,\xi)$ to problem \eqref{eq:DIRICHLESTATIONARYTWS}-\eqref{eq:LIMITCONDTWPROFILES} posed in $\Omega = D\times\RR$ as $j \to +\infty$, with speed of propagation $c_{\ast}$. It depends only on $p$ and $D$ and it is the minimal speed w.r.t. all finite TW solutions to problem \eqref{eq:DIRICHLESTATIONARYTWS}-\eqref{eq:LIMITCONDTWPROFILES}.
\end{lem}
\paragraph{Proof.} Following the ideas of the proof of Theorem \ref{Theorem:WAVECONSTRUCTION} (cfr. with \emph{Step 4}), we normalize the sequence $\varphi_j = \varphi_j(z,\xi)$ (of solutions to problem \eqref{eq:DIRICHLESTATIONARYTWS}-\eqref{eq:LIMITCONDTWPROFILES} posed in $\Omega_j = D_j\times\RR$) by setting $\varphi_j(0,0) = \Phi_{D_j}/2(0)$ for any $j \in \NN$ (here $\Phi_{D_j} = \Phi_{D_j}(z)$ denotes the unique nonnegative weak solution to \eqref{eq:STATIONARYPROBLEM} posed in $D_j$).

\noindent So, if $c_{j\ast}$ is the sequence of speeds corresponding to $\varphi_j$, it is nondecreasing by Lemma \ref{Lemma:MONOTONICITYOFCASTWRTDOMAINS}, and bounded by \emph{Step 4} of \ref{Theorem:WAVECONSTRUCTION}. We thus deduce the existence of a speed $c_{\ast} > 0$ such that $c_{j\ast} \to c_{\ast}$ as $j \to +\infty$, to which corresponds a TW $\varphi = \varphi(z,\xi)$ in $\Omega$. From Lemma \ref{Lemma:MONOTONICITYOFCASTWRTDOMAINS}, we get that $c_{\ast}$ is minimal. Indeed, if $D_1' \subset D_2' \subset \ldots \subset D$ is another sequence of subsets approximating $D$ and $\varphi_j' = \varphi_j'(z,\xi)$ is the corresponding approximating sequence with speeds $c_{j\ast}'$, we can build a new sequence of sets $E_1 \subset E_2 \subset \ldots \subset D$ with $E_j \in \{D_j\}_{j\in\NN}\cup\{D_j'\}_{j\in\NN}$ and apply Lemma \ref{Lemma:MONOTONICITYOFCASTWRTDOMAINS} to $\{E_j\}_{j\in\NN}$, obtaining that $c_{\ast}$ does not depend on the approximating sequence. This conclude the proof of the lemma. $\Box$
\begin{cor}\label{Corollary:APPROXIMATINGCASTFROMOUTSIDE}
Let $D \subset \RR^N$ be a bounded domain and $D_1 \supset D_2 \supset \ldots \supset D$ a sequence of subsets such that
\[
D = \bigcap_{j \in \NN} D_j.
\]
Then there exists a sequence $\varphi_j = \varphi_j(z,\xi)$ of solutions to problem \eqref{eq:DIRICHLESTATIONARYTWS}-\eqref{eq:LIMITCONDTWPROFILES} posed in $\Omega_j = D_j\times\RR$ which converge to a solution $\varphi = \varphi(z,\xi)$ to problem \eqref{eq:DIRICHLESTATIONARYTWS}-\eqref{eq:LIMITCONDTWPROFILES} posed in $\Omega = D\times\RR$ as $j \to +\infty$, with speed of propagation $c^{\ast}$. It depends only on $p$ and $D$ and it is the maximal speed w.r.t. all finite TW solutions to problem \eqref{eq:DIRICHLESTATIONARYTWS}-\eqref{eq:LIMITCONDTWPROFILES}.
\end{cor}
\paragraph{Proof.} The proof is similar to the above one.

\begin{lem}\label{Lemma:Uniqueness}
The speed of any finite TW solution to problem \eqref{eq:DIRICHLESTATIONARYTWS}-\eqref{eq:LIMITCONDTWPROFILES} is unique. Furthermore, for any couple of finite TWs $\varphi_1 = \varphi_1(z,\xi)$ and $\varphi_2 = \varphi_2(z,\xi)$, there exist $l_1,l_2 \in \RR$ such that
\begin{equation}\label{eq:INEQUALITYWAVEPROFILES}
\varphi_1(z,\xi + l_2) \leq \varphi_2(z,\xi) \leq \varphi_1(z,\xi + l_1), \quad \text{ for any } (z,\xi) \in \Omega.
\end{equation}
\end{lem}
\paragraph{Proof.} We follow the proof of Lemma \ref{Lemma:MONOTONICITYOFCASTWRTDOMAINS} in the case in which both $\varphi_1(\cdot)$ and $\varphi_2(\cdot)$ are finite. So, we consider again the functions $u_i = u_i(z,y,t)$ defined in \eqref{eq:TWSORIGINALUI} and their perturbations $u_{1h}(z,y,t) = u_1(z,y + h,t)$ and $u_2^{\theta}(z,y,t) = u_2(z,y,t-\theta)$ for $0 < \theta < 1$ and suitable $h > 0$. To show $c_{1\ast} \leq c_{2\ast}$, it is enough to establish the validity of inequality \eqref{eq:INEQUALITYFORPERTURBATIONSUI}.

\noindent W.r.t. the previous case, we now have $D_1 = D_2 = D$, and thus, even if the comparison at the boundary is trivial, the comparison a $y = -\infty$ is not. To do it, we crucially employ  \eqref{eq:CONVERGENCEINRELATIVERROR}, which, in particular, tells us that for any $\varepsilon > 0$, $\varphi(z,\xi) \geq (1 - \varepsilon)\Phi(z)$ for any $z \in D$ and $\xi \leq \xi_{\varepsilon}$ (and any (finite) TW to \eqref{eq:DIRICHLESTATIONARYTWS}-\eqref{eq:LIMITCONDTWPROFILES}). Consequently, we get
\[
(1 - \theta)^{-\frac{1}{p-2}} \varphi_2(z,y - c_{2\ast}\ln (1 - \theta)) \geq \frac{1-\varepsilon}{(1 - \theta)^{1/(p-2)}} \Phi(z) \geq \Phi(z) \geq \varphi_1(z, y + h),
\]
for any $\varepsilon \leq 1 - (1-\theta)^{\frac{1}{p-2}}$ and any $y \leq -y_{\varepsilon}$ for some $y_{\varepsilon} > 0$ large enough. Finally, taking $h > 0$ large enough we see that \eqref{eq:INEQUALITYFORPERTURBATIONSUI} is satisfied and we conclude $c_{1\ast} \leq c_{2\ast}$ as before. Note that since both waves are finite, we can change the role of $\varphi_1(\cdot)$ and $\varphi_2(\cdot)$ to obtain the reverse inequality and conclude the proof (note that inequality \eqref{eq:INEQUALITYWAVEPROFILES} easily follows by comparison in the limit $t \to +\infty$, we the effect of $0 < \theta < 1$ becomes negligible). $\Box$

%
%
%
%
%
%
%
%
%
%
%
\section{Long time behaviour for general solutions to (\ref{eq:PLETUBULARDOMAIN})}\label{Section:ProofMainTheorem}
Let us consider a solution $u = u(y,z,t)$ to \eqref{eq:PLETUBULARDOMAIN} with nonnegative and nontrivial initial datum $u_0 \in \mathcal{C}_c(\Omega)$, and its re-normalized version $v = v(y,z,\tau)$ defined in \eqref{eq:SECONDRENORMALIZATIONFORMULA}, with $v(x,0) = u_0(x)$ for all $x \in \Omega$.

\noindent Then, from standard results about $p$-Laplacian diffusion, the solution $v = v(y,z,\tau)$ has an expanding-in-time support that, after a waiting time $T =T(\Omega,p,N,u_0) \geq 0$, touches the fixed boundary $\partial \Omega$. From that moment, the \emph{free boundary} is composed by two disjoint sets
\[
S_v^{\pm}(\tau) := \{(z,S^{\pm}(z,\tau)): z \in D \}, \quad \tau > T
\]
where
\[
S_v^+(z,\tau) := \inf\{y > 0: v(z,y,\tau) = 0 \}, \qquad S_v^-(z,\tau) := \sup\{y < 0: v(z,y,\tau) = 0 \},
\]
defined for any $z \in D$ and $\tau > T$, and the following property holds:
\[
v(z,y,\tau) > 0 \quad  \Leftrightarrow  \quad S_v^-(z,\tau) < y < S_v^+(z,\tau), \quad \text{ for any } z \in D, \; \tau > T.
\]
We are ready to prove our main theorem.

\paragraph{Proof.} Let $u_0 \in \mathcal{C}_c(\Omega)$ be nonnegative and nontrivial, and $v = v(z,y,\tau)$ the corresponding solution to \eqref{eq:REACTIONTRANSFORMATION}. The proof of the theorem is based on suitable comparison techniques which employ the TW solution to \eqref{eq:DIRICHLESTATIONARYTWS}-\eqref{eq:LIMITCONDTWPROFILES} in three different ways.

\emph{Step 1: Proof of part (i).} Let us fix $0 < c < c_{\ast}$. For any number $0 < \varepsilon < 1$ we define $D_{\varepsilon} := \delta_{\varepsilon} D$, with $\delta_{\varepsilon} := (1-\varepsilon)^{-\frac{p-2}{p}} > 1$, which implies $D \subset D_{\varepsilon}$. Out of clarity, we postpone  the choice of $0 < \varepsilon < 1$ which, for the moment, can be thought as a fixed parameter in $(0,1)$.

\noindent We want to compare the two nonnegative weak solutions to \eqref{eq:STATIONARYPROBLEM} posed in $D$ and $D_{\varepsilon}$, namely $\Phi = \Phi(z)$ and $\Phi_{\varepsilon} = \Phi_{\varepsilon}(z)$, respectively. To do that, we define
\[
\underline{\Phi}(z) := A_{\varepsilon} \Phi_{\varepsilon}(\delta_{\varepsilon} z), \quad z \in D.
\]
Now, straightforward computations show that $\underline{\Phi} = \underline{\Phi}(z)$ is a sub-solution to problem \eqref{eq:STATIONARYPROBLEM} posed in $D$ (note that the comparison at the boundary $\partial D$ is trivial since $\Phi_{\varepsilon} = 0$ in $\partial D_{\varepsilon}$) if and only if $0 \leq A_{\varepsilon} \leq 1 -\varepsilon$. Consequently, the elliptic comparison principle gives us
\[
\Phi(z) \geq A_{\varepsilon} \Phi_{\varepsilon}(\delta_{\varepsilon}z), \quad z \in D.
\]
Thus, taking $A_{\varepsilon} = 1 -\varepsilon$, we write
\[
1 - \varepsilon = \frac{1-a_{\varepsilon}}{1-b_{\varepsilon}},
\]
for some $a_{\varepsilon},b_{\varepsilon} \in (0,1)$ with $a_{\varepsilon},b_{\varepsilon} \sim 0$ for $\varepsilon \sim 0$ (for instance $a_{\varepsilon} = \varepsilon(3-\varepsilon)/2$ and $b_{\varepsilon} = \varepsilon/2$), so that we get the first fundamental inequality
\begin{equation}\label{eq:FUNDAMENTALINEQSTEADYSTATES}
\left(1 - b_{\varepsilon}\right) \Phi(z) \geq (1 - a_{\varepsilon}) \Phi_{\varepsilon}(\delta_{\varepsilon}z), \quad z \in D.
\end{equation}
Note that from the definition of $a_{\varepsilon}$ and $b_{\varepsilon}$, it automatically follows $\varepsilon < a_{\varepsilon} < 1$ and $0 < b_{\varepsilon} < a_{\varepsilon}$. We recall that, from Lemma \ref{LEMMACONVERGENCEINBOUDEDDOMAINS}, there is $\tau_{\varepsilon} > 0$ which gives us the second main inequality:
\begin{equation}\label{eq:SECONDFUNDAMENTALINEQSTEADYSTATES}
v(z,y,\tau+\tau_{\varepsilon}) \geq \left(1 - b_{\varepsilon}\right)\Phi(z), \quad  z \in D, \; 0 \leq y \leq 1, \; \tau \geq 0.
\end{equation}
Now, let us consider the function
\[
\underline{v}(z,y,\tau) := (1 - a_{\varepsilon})\,\varphi_{\varepsilon}(\delta_{\varepsilon}z,\delta_{\varepsilon}y - c_{\varepsilon}\tau + l), \quad z \in D, \; y,\tau \geq 0, \; l \geq 0,
\]
where $\varphi_{\varepsilon} = \varphi_{\varepsilon}(z,\xi)$ is the wave solution to \eqref{eq:DIRICHLESTATIONARYTWS}-\eqref{eq:LIMITCONDTWPROFILES} posed in $\Omega_{\varepsilon} := D_{\varepsilon} \times (0,\infty)$ and its speed is defined by
\[
c_{\varepsilon} := c_{\ast} \left(\frac{1 - a_{\varepsilon}}{1 - \varepsilon}\right)^{p-2}.
\]
Note that $c_{\varepsilon} < c_{\ast}$ since $\varepsilon < a_{\varepsilon}$, and $c_{\varepsilon} \to c_{\ast}$ as $\varepsilon \to 0$ (this follows from Corollary \ref{Corollary:APPROXIMATINGCASTFROMOUTSIDE}). The last formula allows us to choose $0 < \varepsilon < 1$. Since $a_{\varepsilon} \to 0$ as $\varepsilon \to 0$, we can fix $\varepsilon$ small enough such that $c < c_{\varepsilon} < c_{\ast}$.

\noindent Now, the main fact is that $\underline{v} = \underline{v}(z,y,\tau)$ is a sub-solution to problem \eqref{eq:REACTIONTRANSFORMATION}. Indeed, writing $c_{\ast} = c_{\varepsilon\ast} - o_{\varepsilon}$, where $c_{\varepsilon \ast} > c_{\ast}$ is the speed of the unique wave solution $\varphi_{\varepsilon}$ to \eqref{eq:DIRICHLESTATIONARYTWS}-\eqref{eq:LIMITCONDTWPROFILES} posed in $\Omega_{\varepsilon}$ (the fact that $c_{\varepsilon \ast} > c_{\ast}$ follows by Lemma \ref{Lemma:MONOTONICITYOFCASTWRTDOMAINS}), we have
\[
\begin{aligned}
\partial_{\tau} \underline{v} & = - c_{\ast} \frac{(1 - a_{\varepsilon})^{p-1}}{(1 - \varepsilon)^{p-2}} \partial_{\xi} \varphi_{\varepsilon} = - c_{\varepsilon\ast} \frac{(1 - a_{\varepsilon})^{p-1}}{(1 - \varepsilon)^{p-2}} \partial_{\xi} \varphi_{\varepsilon} +  o_{\varepsilon} \frac{(1 - a_{\varepsilon})^{p-1}}{(1 - \varepsilon)^{p-2}} \partial_{\xi} \varphi_{\varepsilon} \\
& \leq - c_{\varepsilon\ast} \frac{(1 - a_{\varepsilon})^{p-1}}{(1 - \varepsilon)^{p-2}} \partial_{\xi} \varphi_{\varepsilon} = \frac{(1 - a_{\varepsilon})^{p-1}}{(1 - \varepsilon)^{p-2}} \Delta_p\varphi_{\varepsilon} + \frac{(1 - a_{\varepsilon})^{p-1}}{(1 - \varepsilon)^{p-2}} \frac{\varphi_{\varepsilon}}{p-2} \\
& = \Delta_p \underline{v} + \left(\frac{1 - a_{\varepsilon}}{1 - \varepsilon}\right)^{p-2} \frac{\underline{v}}{p-2} \\
& \leq \Delta_p \underline{v} + \frac{\underline{v}}{p-2}, \quad z \in D, \; y,\tau \geq 0,
\end{aligned}
\]
where we used the nonnegativity of $\varphi_{\varepsilon} = \varphi_{\varepsilon}(z,\xi)$, the fact that $\partial_{\xi}\varphi_{\varepsilon} \leq 0$, and the scaling of the $p$-Laplacian, together with our choice of the parameter $\varepsilon < a_{\varepsilon} < 1$.
On the other hand, we have by construction
\begin{equation}\label{eq:THIRDFUNDAMENTALINEQSTEADYSTATES}
(1 - a_{\varepsilon})\Phi_{\varepsilon}(\delta_{\varepsilon} z) \geq \underline{v}(z,y,\tau), \quad z \in D, \; y, \tau \geq 0, \; l \geq 0,
\end{equation}
which is our third fundamental inequality. Consequently, thanks to fact that $\varphi_{\varepsilon} = \varphi_{\varepsilon}(z,\xi)$ has bounded support to the right, we can combine inequalities \eqref{eq:FUNDAMENTALINEQSTEADYSTATES}, \eqref{eq:SECONDFUNDAMENTALINEQSTEADYSTATES}, and \eqref{eq:THIRDFUNDAMENTALINEQSTEADYSTATES}, and properly choose $l = l_{\varepsilon} \geq 0$ such that
\[
v(z,y,\tau_{\varepsilon}) \geq \underline{v}(z,y,0), \quad z \in D, \; y \geq 0.
\]
Moreover, from the same inequalities we easily get $v(z,0,\tau+\tau_{\varepsilon}) \geq \underline{v}(z,0,\tau)$ for any $z \in D$, $\tau \geq 0$ and, from the fact that
\[
v(z,y,\tau+\tau_{\varepsilon}) = 0 = \underline{v}(z,y,\tau), \quad  z \in \partial D, \; y,\tau \geq 0,
\]
we deduce by comparison
\[
v(z,y,\tau + \tau_{\varepsilon}) \geq (1 - a_{\varepsilon})\,\varphi_{\varepsilon}(\delta_{\varepsilon}z,\delta_{\varepsilon}y - c_{\varepsilon}\tau + l_{\varepsilon}), \quad z \in D, \; y,\tau \geq 0.
\]
Since $a_{\varepsilon} \sim 0$ for $\varepsilon \sim 0$, the thesis follows by taking the limit as $\tau \to +\infty$ and recalling that $0 < c < c_{\varepsilon} < c_{\ast}$, together with $\lim_{\xi \to -\infty}\varphi_{\varepsilon}(z,\xi) = \Phi_{\varepsilon}(z) \geq \Phi(z)$ for any $z \in D$.

\emph{Step 2: Proof of part (ii).} In this case we can easily get a bound from above which uses the finite wave $\overline{v}(z,y,\tau) = \varphi(z,y - c_{\ast}\tau - l)$. Indeed, thanks to the universal estimate \eqref{eq:UNIVERSALESTIMATE}, we can assume $u_0(\cdot)$ smaller, for instance $u_0(z,y) \leq \Phi(z)/2$. This allows us to choose $l \geq 0$ so that $\overline{v}(z,y,0) = \varphi(z,y - l) \geq u_0(z,y)$ for any $(z,y) \in \Omega$, and so $\overline{v}(z,y,\tau) = \varphi(z,y - c_{\ast}\tau - l) \geq v(z,y,\tau)$ for any $(z,y) \in \Omega$ and $\tau > 0$ (note that the comparison at the boundary $\partial\Omega$ is trivial). Thus we get the assertion (ii), since each point of the free boundary of $\varphi = \varphi(z,\xi)$ moves with speed $c_{\ast} > 0$.

\noindent We conclude the proof by pointing out that in both part (i) and (ii), we have focused on points belonging to the half-line $y \geq 0$. The same methods apply to the set $y \leq 0$ by reflection $y \to -y$. $\Box$

%
%
%
%
%
%
%
%
%
%
%
\section{Comments and open problems}\label{Section:Finalcomments}
We end the paper with some comments and open problems.
\paragraph{The linear case $\boldsymbol{p = 2}$.} Let us briefly discuss the simpler linear framework, since it presents interesting differences from the nonlinear one $p > 2$ treated here. In the linear  case we should look for re-scaled or re-normalized TWs to the problem
\begin{equation}\label{eq:HETUBULARDOMAIN}
\begin{cases}
\partial_tu = \Delta u \quad &\text{ in } \Omega \times 0,\infty) \\
u = 0 \quad &\text{ in } \partial\Omega \times (0,\infty),
\end{cases}
\end{equation}
where as before $\Omega = D \times \RR$. However, as already observed in \cite{Vazquez2007:art}, we are left without a precise scale of re-normalization, due to the lack of a ``universal estimate'' like \eqref{eq:UNIVERSALESTIMATE}, so that it is not clear if a re-scaled version of \eqref{eq:REACTIONTRANSFORMATION} is available or not. We find two interesting particular cases.

$\bullet$ A natural ansatz is
\[
u(z,y,t) = e^{-\lambda_1 t} \Phi(z) w(y,t),
\]
where $\lambda_1 > 0$ is the first eigenvalue of the Dirichlet Laplacian on $D$, $\Phi = \Phi(z)$ the corresponding eigenfunction and $w = w(y,t)$ satisfies the one-dimensional Heat Equation posed in the whole line
\[
\partial_t w = \partial_{yy} w \quad \text{ in } \RR \times (0,\infty).
\]
Travelling wave solutions to such equation exist and have the form
\[
w(y,t) = A e^{c(ct-y)} + B,
\]
where $A,B \in \RR$ are two free parameters and $c \geq 0$ is the wave speed (note that each speed is admissible). We thus obtain
\[
u(z,y,t) = e^{-\lambda_1 t} \Phi(z) \left[ Ae^{c(ct-y)} + B \right],
\]
which is stationary for $c^2 = \lambda_1$ and $B = 0$, but $u$ is always unbounded in the $y$ variable. Thus, TW solutions for the ``tubular Heat Equation'' exist, but come from an infinity of the l.h.s end of the tube which is not admissible in our problem setting.

$\bullet$ Another possibility, that is closer to our  work for $p>2$, is taking
\[
w(y,t) = \frac{1}{\sqrt{4 \pi t}} \, e^{-\frac{y^2}{4t}},
\]
which gives the re-scaled solution
\[
v(z,y,t) := e^{\lambda_1 t} u(z,y,t) =  \frac{\Phi(z)}{\sqrt{4 \pi t}} \, e^{-\frac{y^2}{4t}}.
\]
This solution decays as $O(t^{-1/2})$ and does not have a long time behaviour of TW type, but a self-similar one instead. Eliminating the $t^{1/2}$ factor, the level sets move with law $y \sim \sqrt{t}$, which is in sharp contrast with the logarithmic one, $y \sim c_{\ast}\ln t$, of the $p$-Laplacian and Porous Medium slow diffusion ones. As mentioned in the introduction, this interesting difference is due to the fact that in the slow diffusion setting the homogeneous Dirichlet conditions play a stronger role in the loss of mass through the lateral boundary.
\paragraph{One-sided propagation.} Through all the paper we have focused on solutions $u = u(z,y,t)$ with nonnegative initial data $u_0 \in \mathcal{C}_c(\Omega)$. Here we show that for a class of initial data which are positive on the left and compactly supported on the right a stronger bound on the corresponding solution holds and the asymptotic behaviour is much different. As always, we state our result for solutions to the more natural problem \eqref{eq:REACTIONTRANSFORMATION}, whose proof follows from an easy comparison.
\begin{thm}(One-sided propagation)
Let $v_0 = v_0(z,y)$ be a nonnegative and continuous function satisfying
\[
\varphi(z,y + l_1) \leq v_0(z,y) \leq \varphi(z,y + l_2),
\]
for some $l_1,l_2 \in \RR$, where $\varphi = \varphi(z,\xi)$ is the TW found in Theorem \ref{ASYMPTOTICBEHAVIOURTHEOREM}. Then the solution $v = v(z,y,\tau)$ to \eqref{eq:REACTIONTRANSFORMATION} with initial data $v_0$ satisfies
\[
\varphi(z,y - c_{\ast}\tau + l_1) \leq v(z,y,\tau) \leq \varphi(z,y - c_{\ast}\tau + l_2),
\]
where $c_{\ast} > 0$ is the critical speed corresponding to $\varphi = \varphi(z,\xi)$.
\end{thm}
\normalcolor
\paragraph{The fast diffusion range.} Problem \eqref{eq:PLETUBULARDOMAIN} can be posed in the fast diffusion range $1<p < 2$, i.e., in the singular diffusion framework. It is well-known (cfr. for instance with \cite{V1:book}) that in this range the behaviour of solutions strongly differs from the case $p > 2$. One of the main facts is that solutions spread through the space with infinite speed of propagation and, consequently, they do not have a free boundary. In a recent work \cite{AA-JLV2:art}, we studied the Fisher-KPP problem with $p$-Laplacian diffusion in the ``good range'' of fast diffusion $2N/(N+1) < p < 2$ and we proved that there are not TW solutions, but general solutions propagate exponentially fast for large times, while if $1 < p < 2N/(N+1)$ solutions vanish in finite time (whether the initial datum is nonnegative and compactly supported). When $\Omega = D \times \RR$ is a tube, things seem to work in a quite different way and solutions may extinguish in finite time for any $1 < p < 2$. This peculiar feature comes from the fact that problem \eqref{eq:PLEBOUNDEDDOMAIN} admits solutions with extinction in finite time: $u(\cdot,t) \to 0$ uniformly as $t \to T$, for some finite time $T > 0$ (cfr. with \cite{Bon-VazPME2010:art,Bon-Grillo-Vaz:art} for the fast Porous Medium framework and \cite{Bon-VazPLE2010:art} for the fast $p$-Laplacian one) that could be employed as super-solutions to \eqref{eq:PLETUBULARDOMAIN}, obtaining that also solutions  defined in tubes extinguish in finite time. Clearly, in this case the study of propagation of solutions has a very different nature. \normalcolor
\paragraph{Sharp asymptotics.} Theorem \ref{ASYMPTOTICBEHAVIOURTHEOREM} shows that solutions to problem \eqref{eq:REACTIONTRANSFORMATION} propagate along the tube with constant speed of propagation $c_{\ast} = c_{\ast}(p,D) > 0$ for large times. However, it does not contain any precise information on the limit (if it exists) of general solutions apart from the level set information. In other words, given a solution $v = v(z,y,\tau)$ to \eqref{eq:REACTIONTRANSFORMATION} with nonnegative initial data $v_0 = v_0(z,y)$, it is not clear if there is a limit profile $v_{\infty} = v_{\infty}(z,y,\tau)$ such that $v \to v_{\infty}$ as $\tau \to +\infty$, where the convergence is intended in some suitable sense.

A reasonable guess is that the limit exists and  $v_{\infty}(z,y,\tau) = \varphi(z,y - c_{\ast}\tau + l)$, where $\varphi = \varphi(z,\xi)$ is the wave solution corresponding to $c_{\ast}$ and $l \in \RR$ is a suitable shift depending on the data. This problem is much studied in the context of reaction-diffusion equations and the existence of a limit and its properties strongly depend on the initial data (cfr. for instance with \cite{B2:art,FifeMcLeod1977:art,Hamel-N-R-R:art} for the linear diffusion framework). In the nonlinear setting much less has been done (we quote the work \cite{Du-Quiros-Zhou:art} for the Fisher-KPP setting and Porous Medium diffusion and \cite{Garriz2018:art} for the bistable counterpart). The $p$-Laplacian case is completely open. Furthermore, we stress that the problem posed in tubular domains could be of particular interest since it intrinsically possesses features of both $1$-dimensional and $N$-dimensional problems. As showed in \cite{Du-Quiros-Zhou:art}, the asymptotic behaviour of solutions to the Fisher-KPP problem with Porous Medium diffusion strongly depends on the dimension of the space on which such solutions are defined. When the problem is posed in tubes it is thus not clear what to expect, and the study of the existence of a limit seems to be an interesting open problem.

\paragraph{The problem of a free boundary sliding on a wall.}
A number of degenerate diffusion equations exhibit the property of Finite Propagation
 whereby solutions with initial data localized in a region may expand the
support in time, but only to a finite distance of the original support for every finite time. This happens in the present situation, since the nonnegative solutions of the $p$-Laplacian equation in the tube
$\Omega = D \times \RR$, with $p>2$ and, say, bounded and compactly supported initial  data, evolve in time so that the support spreads to eventually reach every point of the domain. Moreover, the support reaches the lateral boundary $\partial \Omega = \partial D\times \RR$ in finite time, and then spreads along the tube. This situation reminds  of a main problem  in  fluid mechanics which concerns the way fluids slide along confining walls.

Since we have constructed a special solution (of the logarithmic-time travelling wave type), and this solution exhibits a finite free boundary, the problem is then to determine the geometry and regularity of this free boundary and, more precisely, the way it makes contact with the lateral boundary $\partial\Omega$. Numerical and formal calculations show that in the case of zero Dirichlet boundary conditions such a contact
must be tangential, thus  eliminating the existence of a contact angle. We need to rigorously prove that conjecture  and to understand the typical behaviour at the contact points. This is a difficult open problem.

The mentioned problem of the flow front sliding along a wall is equally posed for other equations and systems. A main example is  the Porous Medium Equation, $u_t = \Delta u^m$ with $m > 1$, where the same open problem was posed in \cite{Vazquez2004:art}. It has not been solved.

\vskip .5cm


\noindent {\textbf{\large \sc Acknowledgments.}} The first author has been partially funded by the ERC Advanced Grant 2013 n.~339958 ``Complex Patterns for Strongly Interacting Dynamical Systems - COMPAT'' and the GNAMPA project ``Ottimizzazione Geometrica e Spettrale'' (Italy). The second author is partially funded by Project MTM2014-52240-P (Spain) and is an Honorary Professor at Univ. Complutense de Madrid, where part of this work was done. We thank David G\'omez Castro, of UCM, for the numerical computations involved in Figure 1.

\


\noindent 2010 \textit{Mathematics Subject Classification.}
35K57,  
35K65, 
35C07,  	
35K55, 

\medskip

\noindent \textit{Keywords and phrases.} Travelling waves, $p$-Laplacian diffusion in tubes, Long-time behaviour, Re-normalized variables.

\

%
%
%

%
%
%

%

\

\

\end{document}